\documentclass{amsart}
\usepackage{amscd}
\usepackage{amssymb}
\theoremstyle{plain}
\newtheorem{theorem}{Theorem}[section]

\newtheorem{proposition}[theorem]{Proposition}
\newtheorem{corollary}[theorem]{Corollary}
\newtheorem{algorithm}[theorem]{Algorithm}

\theoremstyle{definition}
\newtheorem{definition}[theorem]{Definition}
\newtheorem{def-prop}[theorem]{Definition-Proposition}
\newtheorem{example}[theorem]{Example}

\theoremstyle{remark}
\newtheorem{remark}{Remark}[section]

\begin{document}
\title[Locally semi-simple representations of quivers]
{ Locally semi-simple representations of quivers.}
\address{117437, Ostrovitianova, 9-4-187, Moscow, Russia.}
\email{mitia@mccme.ru}
\author{D.A.Shmelkin}
\begin{abstract}
We suggest a geometrical approach to 
the semi-invariants of quivers based on Luna's
slice theorem and the Luna-Richardson theorem.
The locally semi-simple representations are defined
in this spirit but turn out to be connected with 
stable representations in the sense of GIT, 
Schofield's perpendicular categories,
and Ringel's regular representations.
As an application of this method we obtain an independent short proof 
of the theorem of Skowronsky and Weyman about semi-invariants
of the tame quivers.
\end{abstract}


\maketitle
\section{Introduction}
Let $Q$ be a finite quiver, i.e., an oriented graph. 
We fix the notation as follows: denote by $Q_0$ and $Q_1$ the sets
of the vertices and the arrows of $Q$, respectively. For any arrow
$\varphi\in Q_1$ denote by $t\varphi$ and $h\varphi$ its tail and
its head, respectively. A representation $V$ of $Q$ over 
an algebraically closed field ${\bf k}$, ${\rm char}\,{\bf k}=0$,
consists in defining  a vector space $V(i)$ over ${\bf k}$ 
for any $i\in Q_0$ and  a ${\bf k}$-linear map 
$V(\varphi):V(t\varphi)\to V(h\varphi)$ for any 
$\varphi\in Q_1$. The dimension vector $\dim V$ is the collection
of $\dim V(i), i\in Q_0$. For a fixed dimension $\alpha$ we may
set $V(i)={\bf k}^{\alpha_i}$. Then the set $R(Q,\alpha)$ of 
the representations of dimension $\alpha$ is converted into the
vector space
\begin{equation}
R(Q,\alpha) = \bigoplus_{\varphi\in Q_1} {\rm Hom}
({\bf k}^{\alpha_{t\varphi}},{\bf k}^{\alpha_{h\varphi}}).
\end{equation}
A homomorphism $H$ of a representation $U$ of $Q$ to another representation,
$V$ is a collection of linear maps $H(i), U(i)\to V(i)\in Q_0$ such that
for any $\varphi\in Q_1$ holds $V(\varphi) H(t\varphi) = H(h\varphi)U(\varphi)$.
The endomorphisms, automorphisms, and isomorphisms are defined naturally.
Hence, the isomorphism classes
of representations of $Q$ are the orbits of a reductive group
$GL(\alpha) = \prod_{i\in Q_0} GL(\alpha_i)$ acting naturally on 
$R(Q,\alpha)$: 
$(g(V))(\varphi) = g(h\varphi) V(\varphi) (g(t\varphi))^{-1}$.
Set $SL(\alpha)=\prod_{i\in Q_0} SL(\alpha_i)\subseteq GL(\alpha)$.

Assume that $Q$ has no oriented cycles. Then for any dimension 
$\alpha$ the algebra ${\bf k}[R(Q,\alpha)]^{GL(\alpha)}$ of 
$GL(\alpha)$-invariant regular functions is trivial and
the unique $GL(\alpha)$-closed orbit is the origin of
$R(Q,\alpha)$. It is however interesting to study
the $SL(\alpha)$-invariant functions or the semi-invariants
of $GL(\alpha)$. The aim of this paper is to suggest
a geometrical approach to this study in the spirit of Luna's
papers \cite{lu1},\cite{lu2}. For this we need to describe
the closed orbits of $SL(\alpha)$. Consider a more general 
setting of a connected reductive group $G$ acting on an affine variety $X$,
$x\in X$, $G'\subseteq G$ is the commutant.
In \ref{lss-criteria} we prove that $G'x$ is closed in $X$ if and only if 
$Gx$ is closed in an open affine neighborhood $X_f\subseteq X$, where
$f$ is semi-invariant. 
We call such $x$ a {\it locally semi-simple} point. 
In \ref{gen_stab_lss} we prove that there exists a {\it generic stabilizer
of locally semi-simple points} and in \ref{lr_lss} 
we obtain a special version  of the Luna-Richardson theorem (see \cite{lu2})
that can be called the Luna-Richardson theorem about semi-invariants.

Locally semi-simple representations of quivers
turn out to be closely connected with the 
stable representations in the sense of GIT (see \cite{mf},\cite{ki})
and the perpendicular categories
introduced in \cite{sch}. Namely we prove in \ref{crit_lss}
that a representation is locally semi-simple if and only if
it is a sum of simple objects in a perpendicular category;
actually this is just a more strong version of \cite[Proposition~3.2]{ki}.

Recall (\cite{kac}) that a decomposition 
$\alpha = \sum_{i=1}^t \beta_i$ with 
$\alpha,\beta_1,\cdots,\beta_t$ $\in {\bf Z}_+^{Q_0}$ is called
{\it canonical} if $\beta_1,\cdots,\beta_t$ are Schur roots and 
the set of representations $R_1+\cdots+R_t$ such that $R_i$ is indecomposable 
and $\dim R_i = \beta_i$ contains an open dense subset in $R(Q,\alpha)$. 
On the other hand, Ringel applied in \cite{ri} the term "canonical decomposition" 
for a different notion and we need that notion too.  
We therefore call the decomposition introduced by Kac {\it generic}
(as e.g. in \cite{skw}).

It is well-known that the generic decomposition
corresponds to the generic stabilizer in the sense that
the torus $T\subseteq GL(\alpha)$ of rank
$t$ naturally corresponding to this decomposition is a maximal torus
in a generic stabilizer for the action of $GL(\alpha)$.
Analogically, the maximal torus of the generic stabilizer
of locally semi-simple points yields another decomposition of
$\alpha$ that we call {\it generic locally semi-simple}.
Given a non-trivial locally semi-simple point $V\in R(Q,\alpha)$, 
we show in \ref{slice_decomp} how to describe both decompositions 
in terms of those for a quiver  $\Sigma_V$ and a dimension vector
$\gamma$. The corresponding linear group  $(GL(\gamma), R(\Sigma_V,\gamma))$
is nothing else but the image of ${\rm Aut}(V)$ under
the {\it slice representation}. Note that for
$V$ {\it semi-simple} such a form of the slice representation 
is known after\cite{lbp}.
Applying  \ref{lr_lss} for the quiver setting, we get
a useful general description \ref{lr_quiver} of 
${\bf k}[R(Q,\alpha)]^{SL(\alpha)}$.

We apply the above methods to the case when $Q$ is a tame quiver.
Ringel introduced in \cite{ri} the category  ${\mathcal R}$ 
of {\it regular} representations.
For the tame quivers he described explicitly the simple objects
of ${\mathcal R}$. 
We prove in \ref{regular_is_perp} that ${\mathcal R}$  is 
the union of the perpendicular categories $^{\perp}S$ with $S$ running
over the homogeneous simple regular representations.
For a dimension vector $\alpha$ such that $R(Q,\alpha)$ contains
a regular representation Ringel introduced a {\it canonical} decomposition
of $\alpha$ into a sum of the dimensions of the simple regular representations.
An important observation \ref{canon_yield_lss} is that
the regular representations $V\in R(Q,\alpha)$ corresponding to this decomposition
are locally semi-simple. Moreover, the linear group  $(GL(\gamma), R(\Sigma_V,\gamma))$
corresponds to a quiver being a disjoint union of equioriented
$A_n$-type quivers. For the latter quivers we can easily describe
the generic and the generic locally semi-simple decomposition.
Together with \ref{slice_decomp} this yields a way to determine both
decompositions for $Q$ and $\alpha$ (in what concerns the generic decomposition
we recover \cite[Theorem~3.5]{ri}).

The algebras of semi-invariants 
of tame quivers $Q$ have been studied in several papers
including \cite{ri},\cite{hh},\cite{schw}. In \cite{skw} Skowronsky
and Weyman proved that ${\bf k}[R(Q,\alpha)]^{SL(\alpha)}$ is a 
{\it complete intersection} for any $\alpha$.
Using \cite{dw} and \ref{lr_quiver}, we give an independent short proof 
of this result (\ref{homog_gen},\ref{main}).

\section{Locally semi-simple points.}\label{sect_lss}
Throughout this section let $G$ denote a connected
reductive group, $G'$ stands for its commutant, $T=G/G'$
is a torus; let $X$ be an irreducible affine variety
acted upon by $G$. For a $T$-character $\chi\in \Xi (T)$
denote by ${\bf k}[X]^{(G)}_{\chi}$ the module
of eigenvectors of $G$ with weight $\chi$, denote
by $G_{\chi}$ the kernel of $\chi$. 
Clearly, the algebra ${\bf k}[X]^{G'}$ of $G'$-invariant
functions on $X$ is the direct sum:
${\bf k}[X]^{G'}=\oplus_{\chi\in \Xi (T)} {\bf k}[X]^{(G)}_{\chi}$. 

\begin{theorem}\label{lss-criteria} The following properties of $x\in X$ are equivalent:

(i) for a semi-invariant $f\in {\bf k}[X]^{(G)}_\chi$,
$f(x)\neq 0$ and $Gx$ is closed in $X_f$

(ii) for a character $\chi\in \Xi (T)$, the orbit $G_{\chi}x$ 
is closed in $X$  

(iii) the orbit $G'x$ is closed in $X$  

(iv) the closure of the orbit $G'x$ in $X$ is contained in $Gx$.
\end{theorem}

\begin{proof}
First observe that if $f\in {\bf k}[X]^{(G)}_\chi$,
$f(x)\neq 0$, and $g\in G$, then 
$g(G_{\chi}x) = Gx \cap \{y\in X\vert f(y) = f(gx)\}$.
This yields the implication $(i)\Rightarrow (ii)$.
Also we note that $Gx$ is a disjoint union of $t(G_{\chi}x)$
with $t$ running over the 1-dimensional coset space $G/G_{\chi}$.
The implication 
$(ii)\Rightarrow (iii)$ follows from the fact that the subgroup
$G'$ is normal in $G_{\chi}$.

Let us prove the implication $(iii)\Rightarrow (i)$.
Let $\pi_{G'}:X\to X/\!\!/G'$ denote the quotient
map. The torus $T$ acts on the quotient $X/\!\!/G'$;
consider a $T$-equivariant embedding of $X/\!\!/G'$ to a $T$-module $W$.
Clearly, if $y=\pi_{G'}(x)$ belongs to $W^T$, then
$Ty=y$ implies $Gx=G'x$ and $(i)$ holds with $f$ being a constant
function. Otherwise,
$y$ is a sum of non-zero $T$-eigenvectors; let $f$ be
the product of the corresponding linear $T$-eigenfunctions.
Then $f$ is a $T$-semi-invariant function on $W$ with respect
to a character $\chi\in \Xi (T)$, and $f(y)\neq 0$.
Moreover, the orbit of $y$ with respect to the kernel
$T_{\chi}$ of $\chi$ is closed, because $y$ is a sum of
$T_{\chi}$-eigenvectors such that the sum of
their characters (with respect to $T_{\chi}$) is zero.  
Consequently, the orbit $G_{\chi}x=T_{\chi}G'x$ is closed,
because is equal to the intersection of the closed pull-back
$\pi_{G'}^{-1}T_{\chi}y$ with the closed subset 
$X_{\dim G'x}=\{z\in X\vert \dim G'z \leq \dim G'x\}$.
So we got $(ii)$ and besides, $f(x)\neq 0$ for the above
$f$ thought of as a $G$-semi-invariant function on $X$.
Set $m = \dim G_{\chi}x$. Then the closure $\overline{Gx}$
of $Gx$ in $X$ is an irreducible variety of dimension $m+1$
and for any $t\in G/G_{\chi}$ the set
$\{y\in \overline{Gx}\vert f(y) = f(tx)\}$ is an
equidimensional closed subvariety in $X$ of dimension $m$
with one of irreducible components being equal to $t(G_{\chi}x)$.
Since $G/G_{\chi}$ acts transitively on the fibers of $f$,
we get $\overline{Gx}_f$ is an equidimensional variety with
one of irreducible components being equal to $Gx$.
However, the dimension of $\overline{Gx}\setminus Gx$ is less than
$\dim Gx$, hence, $\overline{Gx}_f  = Gx$, that is, $Gx$ is closed
in $X_f$.

The implication $(iii)\Rightarrow (iv)$ is obvious.
If $G'x$ is not closed, then  there is $z\in \overline{G'x}$ 
such that $\dim G'z < \dim G'x$. However, for any $z\in Gx$
we have $\dim G'z = \dim G'x$. Therefore $(iv)$ implies $(iii)$.
\end{proof}

The points $x\in X$ such that $Gx$ is closed in $X$
are called {\it semi-simple} in \cite{gv}.
If $X$ is a variety of representations of associative algebras,
then this is not just a definition, since it is proved
that the module corresponding to $x$ is semi-simple if and only if
$Gx$ is closed (see e.g. \cite{kr}). This motivates

\begin{definition}
We call $x\in X$ {\it locally semi-simple} if $x$ fulfills the
equivalent conditions of the above theorem.
\end{definition}

\begin{remark}\label{mum} The property of local semi-simplicity
is intermediate between those of {\it stability} and 
{\it semi-stability} introduced by Mumford (\cite{mf}).
Recall that in our context $x$ is called $\chi$-semistable if
$f(x)\neq 0$ for a non-constant semi-invariant $f\in {\bf k}[X]^{(G)}_\chi$,
and $x$ is called $\chi$-stable if $x$ is $\chi$-semistable, 
the stabilizer of $x$ is equal to the kernel of the action $G:X$, and
the orbit $Gx$ is closed in $X_f$. So locally semi-simple
points meeting condition {\it (i)} of \ref{lss-criteria}
are $\chi$-semistable and $x$ is $\chi$-stable if and only if
$x$ is locally semi-simple with trivial stabilizer.
\end{remark}

For a reductive group $M$ acting on an affine variety $Y$
D.Luna introduced in \cite{lu1},\cite{lu2} the concept
of {\it \'etale slice} at a semi-simple point $y\in Y$.
First of all by Matsushima's criterion \cite{ma}, the stabilizer
$M_y$ is reductive. The Luna slice theorem (\cite{lu1})
states that there exists an  \'etale slice
$S\subseteq Y$ at $y$ such that $S\ni y$ is affine, locally closed,
$M_y$-stable, and the natural map 
$\varphi_y:M*_{M_y}S\to Y, [m,y]\to my$ is {\it excellent}
(see precise definition in \cite{lu1}), in particular the image
of $\varphi_y$ is affine and the restriction of $\varphi_y$
to any fiber of the $M$-quotient map is an isomorphism. 
Further, assume that $Y=V$ is a vector space and $M$ acts on $V$ 
by a linear representation, $v=y$; choose a $M_v$-stable complementary
subspace $N$ to ${\rm T}_v Mv$ in $V$. Then as $S$ we can take
$S = v+ N_0$ for an open affine subset $N_0\subseteq N$ containing 0.
The representation $\sigma_v: M_v\to GL(N)$ is called in this case the 
{\it slice representation} of $v$ and can be calculated by the formula
(${\rm Ad}$ stands for the adjoint representation):
\begin{equation}\label{slice_formula}
\sigma_v = (V + {\rm Ad}M_v)/({\rm Ad}M)\vert_{M_v},
\end{equation}

Let $\pi_{G'}:X\to X/\!\!/G' = {\rm Spec}{\bf k}[X]^{G'}$
denote the quotient map. For $\xi\in X/\!\!/G'$
denote by ${\mathcal O}_{\xi}$ the unique $G'$-closed orbit
in $\pi_{G'}^{-1}(\xi)$. The quotient $X/\!\!/G'$
carries the Luna stratification (\cite{lu1}) by the disjoint
locally closed subvarieties 
$(X/\!\!/G')_{(L)}=\{\xi\in X/\!\!/G'\vert 
{\mathcal O}_{\xi}\cong G'/L\}$, where $L$ is 
a subgroup in $G'$.
We consider a similar stratification with respect to the action of 
$G$, as follows.
For a (reductive) subgroup $M\subseteq G$
denote by $(X/\!\!/G')^G_{(M)}$ the set of all $\xi\in X/\!\!/G'$
such that $G_z$ is $G$-conjugate to $M$ for 
$z\in {\mathcal O}_{\xi}$. Clearly, if the
subgroups $M_1$ and $M_2$ are $G$-conjugate,
then $M_1\cap G'$ and $M_2\cap G'$ are $G'$-conjugate,
hence each Luna stratum
$(X/\!\!/G')_{(L)}$ is a union of $(X/\!\!/G')^G_{(M)}$
with $M\cap G'$ being $G'$-conjugate to $L$.

\begin{proposition}\label{strat_loc_closed}
 $(X/\!\!/G')^G_{(M)}$ is locally closed.
\end{proposition}
\begin{proof}
Apply the slice theorem for $z\in {\mathcal O}_{\xi}^M$ and $G'$.
Since $z$ is $M$-invariant and $M$ normalizes $G'_z$, the slice
$S$ can be chosen to be $M$-stable. Then the map $\varphi_z$
is $M$-equivariant. Denote by 
$\varphi_z/\!\!/G':S/\!\!/G'_z\to X/\!\!/G'$ 
the \'etale covering of a neighborhood of $\xi$ in $X/\!\!/G'$
given by the slice theorem. Then the $M$-stratum is covered
by $S^M$, hence is locally closed.
\end{proof}

\begin{proposition}\label{strat_fin}
The stratification $X/\!\!/G'=\sqcup_{M}(X/\!\!/G')^G_{(M)}$ is finite.
\end{proposition}

\begin{proof}
Since the Luna stratification is finite,
it is sufficient to show that each Luna stratum $(X/\!\!/G')_{(L)}$
is decomposed into finitely many strata $(X/\!\!/G')^G_{(M)}$.
Clearly, we may assume $M\cap G' = L$.
Take $\xi\in (X/\!\!/G')_{(L)}$; then $\xi\in (X/\!\!/G')^G_{(M)}$
with $M/L\cong (G/G')_{\xi}$ so $M$ is an extension
of $L$ by a diagonalizable subgroup in the normalizer
$N_G(L)$. Choose a maximal torus $A\subseteq N_G(L)$.
Then each point $z$ with $G'_z=L$ is 
$N_G(L)$-conjugate to a point $w$ such that the identity component of
$G_w$ is contained in $LA$. On the other hand,
$X^L$ can be divided into finitely many subsets with
constant stabilizer with respect to $A$.
So there are finitely many $N_G(L)$-conjugacy classes
of stabilizers in $G$ of points $z$ with $G'_z=L$.
\end{proof}

Recall that a subgroup $H_0\subseteq G'$ is called
{\it principal isotropy group} if  $(X/\!\!/G')_{(H_0)}$
is the unique open and dense Luna stratum. 
By Propositions \ref{strat_loc_closed} and \ref{strat_fin}
we have:

\begin{def-prop}\label{gen_stab_lss}
A subgroup $H\subseteq G$ is called {\it generic stabilizer of a 
locally semi-simple point} if $(X/\!\!/G')^G_{(H)}$
is the unique open and dense $G$-stratum of $X/\!\!/G'$.
The intersection $H\cap G'$ is a principal isotropy group
and the image of $H$ in $G/G'$ is the kernel of the action $G/G':X/\!\!/G'$.
\end{def-prop}

We now want to describe locally semi-simple points and their
stabilizers in terms of the Luna slice theorem with respect
to the group $G$. Indeed, if $x$ is locally semi-simple and 
$Gx$ is closed in $X_f$ for a semi-invariant $f$, then there is 
an \'etale slice $S$ at $x$ with respect
to $X_f$, and if $X=V$ is a vector space, 
then an \'etale slice of type $x+N_0$ exists. 

\begin{proposition}\label{lss_from_slice}
If $G:V$ is a linear representation, $v\in V$ is a locally semi-simple point,
then for any $n\in N_0$, $v+n$ is locally semi-simple with respect to $G$
if and only if $n\in N$ is locally semi-simple with respect to $G_v$. 
\end{proposition}

\begin{proof}
Assume that $G(v+n)$
is closed in $V_f$ for a $G$-semi-invariant $f$.
Then $\varphi_v^{-1} G(v+n)$ is closed in
$\varphi_v^{-1} V_f = G *_{G_v} (v+(N_0)_{f'})$, where $f'\in {\bf k}[N]$
is defined as $f'(n') = f(v+n'),n'\in N$
so that $f'$ is $G_v$-semi-invariant.
Moreover, $\varphi_v$ is excellent implies that 
$\varphi_v^{-1} G(v+n)$ is a union of finitely many orbits, hence
$G[e,v+n]=G*_{G_v}(v+G_vn)$ is also closed in $G *_{G_v} (v+(N_0)_{f'})$,
equivalently $G_vn$ is closed in $(N_0)_{f'}$.
Since $N_0$ is affine and $N$ is a vector space,
$N_0 = N_d$ for some $d\in {\bf k}[N]$ and $N_0$ is $G_v$-stable
implies that $d$ is $G_v$-semi-invariant. So
$G_vn$ is closed in $N_{f'd}$ and we proved the "only if" part.

Assume that $G_vn$ is closed in $N_{f'}$. As above, we may additionally
assume that $N_{f'}$ is contained in $N_0$. Then by the properties
of an excellent map we have that $G(v+n)$ is closed in an open
subset $V_1$ in the image $V_0$ of $\varphi_v$, such that
$V_0\setminus V_1$ is an equidimensional subvariety of codimension
1 in $V_0$. Since $V_0$ is affine and $V$ is a vector space,
we get $V_1 = V_f$  for some $f\in {\bf k}[V]$, and $V_0$ is
$G$-stable implies $f$ is $G$-semi-invariant.
\end{proof}

It is well-known that the isotropy group $G_v$ for a semi-simple 
$v\in V$ is principal if and only if 
the only semi-simple point in $(G_v,N/N^{G_v})$ is 0.

\begin{corollary}\label{criter_generic}
Let $N = N^{G_v}\oplus N_+$ be a $G_v$-stable decomposition.
Then $v$ is generic if and only if 
the only $G_v$-locally semi-simple point in $N_+$ is 0.
\end{corollary}  

\begin{proof}
If $N_+\setminus\{0\}$ contains a $G_v$-locally semi-simple point,
then $N$ contains a $G_v$-locally semi-simple point $n$
with a proper isotropy subgroup $(G_v)_n\subseteq G_v$. 
Multiplying $n$ by a scalar, we may assume $n\in N_0$, hence
by the Proposition, $v+n$ is $G$-locally semi-simple with stabilizer
$G_{v+n} = (G_v)_n$. Then the closure of 
$(X/\!\!/G')^G_{(G_{v+n})}$ contains $(X/\!\!/G')^G_{(G_v)}$ so 
the closure of the latter can not be equal to $X/\!\!/G'$ and
$v$ is not generic. Conversely, 
if the only $G_v$-locally 
semi-simple point in $N_+$ is 0, then $v$ is a generic locally semi-simple
point in the image $V_0$ of $\varphi_v$. 
Let $V_{lss}\subseteq V$ be the subset of $G$-locally semi-simple points.
By Theorem \ref{lss-criteria}, $V_{lss}$ is also 
the union of $G'$-closed orbits; since $V/\!\!/G'$ is irreducible,
the closure $\overline{V_{lss}}$ also is. Therefore  $\overline{V_{lss}}$ is
the closure of its intersection with $V_0$ and $v$ is generic
in $V$.
\end{proof}

\begin{corollary}\label{generic_lss_from_slice}
If $n\in N_0$ is a generic locally semi-simple point
for the action of $G_v$, then $v+n$ is generic for
$G$ acting on $V$.
\end{corollary}

\begin{proof}
By the slice theorem $G_{v+n} = (G_v)_n$; by the Proposition, 
$v+n$ is locally semi-simple.
Applying formula (\ref{slice_formula}), we get:
$\sigma_{v+n}=\sigma_n$.
Applying Corollary \ref{criter_generic} we conclude the proof.
\end{proof}

The notion of the locally semi-simple point can be used
in order to describe the semi-invariants of $G$. 
Recall that the Luna-Richardson theorem
\cite{lu2} says that if $H_1$ is a principal isotropy group 
for the action $G:X$, then the embedding $X^{H_1}\subseteq X$
gives rise to an isomorphism: 
${\bf k}[X]^G \xrightarrow{\sim}{\bf k}[X^{H_1}]^{N_G(H_1)}$.
Of course one could take $G'$ instead of $G$ and apply the Luna-Richardson
theorem. But there is another way:

\begin{proposition}\label{lr_lss}
Assume that $H$ is a generic stabilizer of a locally semi-simple point
for the action $G:X$. Then the embedding $X^H\subseteq X$
and the group homomorphism $\theta_H:N_{G}(H)/N_{G'}(H)\to G/G'$
give rise to an isomorphism: 
\begin{equation}
{\bf k}[X]^{G'} \xrightarrow{\sim}{\bf k}[X^H]^{N_{G'}(H)}
\end{equation}
of a $\Xi(G/G')$-graded algebra onto a $\Xi(N_{G}(H)/N_{G'}(H))$-graded
algebra. Moreover generic orbit of $N_{G'}(H)$ is closed in $X^H$.
\end{proposition}

\begin{proof}
First prove that we have an isomorphism of algebras.
Since $H$ normalizes $G'$, we have $HG'$ is a reductive subgroup
in $G$. Moreover $H$ is a principal isotropy group for $HG'$ acting
on $X$ so ${\bf k}[X]^{HG'}$ restricts isomorphically
onto ${\bf k}[X^H]^{N_{HG'}(H)}={\bf k}[X^H]^{HN_{G'}(H)} = 
{\bf k}[X^H]^{N_{G'}(H)}$. Also we have: ${\bf k}[X]^{G'}= {\bf k}[X]^{HG'}$.
Indeed for any $h\in H, f\in {\bf k}[X]^{G'}$ we have 
$hf\equiv f$ on $X^H$, hence on $G'X^H$. But $X^H$ intersects
all closed $HG'$-orbits, so any $HG'$-invariant function
is completely defined by its restriction to $HG'X^H=G'X^H$.
Thus we got that restricting $G'$-invariant functions to
$X^H$ we have an isomorphism: 
${\bf k}[X]^{G'} \xrightarrow{\sim}{\bf k}[X^H]^{N_{G'}(H)}$.
Clearly a $G$-semi-invariant function of weight $\chi\in \Xi(G/G')$
restricts to a $N_G(H)$-semi-invariant function of weight
$\theta_H^*(\chi)$. Finally observe that generic point $x\in X^H$
has a closed $HG'$-orbit, because $H$ is a principal isotropy
group of $HG'$. By \cite{lu2}, $N_{HG'}(H)x=N_{G'}(H)x$ is also closed 
in $X^H$. 
\end{proof}
\begin{remark}
The description of ${\bf k}[X]^{G'}$ given by this theorem 
is in general different from the given by the Luna-Richardson
theorem for $G'$. For instance, take as $G'$ the group $SL_2$
acting naturally on ${\bf k}^2\oplus{\bf k}^2$;
as $G$ take the extension of $G'$ by 
$T=\{diag(t,t,s,s)\vert t,s\in{\bf k}^*\}$.  
Then principal isotropy group of $G'$ is trivial
but $H=\{diag(u,1,1,u^{-1})\vert u\in{\bf k}^*\}$.
\end{remark}

\section{Locally semi-simple representations of quivers.}

\begin{definition}
A representation $V$ of a quiver $Q$ is called 
locally semi-simple if $V$ is a locally semi-simple
point of $R(Q,\dim V)$ with respect to $GL(\dim V)$.
\end{definition}

We start with observations as follows:

\begin{proposition}\label{nec_lss}
Assume that $V$ is a locally semi-simple representation.

{\bf 1.} If $V= V_1+V_2$, then both $V_1$ and $V_2$ are locally semi-simple.

{\bf 2.} If $V$ is indecomposable, then ${\rm Aut}(V)={\bf k}^*$.

{\bf 3.} If $V= V_1+V_2$, both $V_1$ and $V_2$ are indecomposable, 
and $V_1\not\cong V_2$, then
we have:
${\rm Hom}(V_1,V_2)=0, {\rm Hom}(V_2,V_1)=0$.
\end{proposition}

\begin{proof}
The assertions 1-3 follow from Theorem \ref{crit_lss} below. 
We give however an independent proof. 
Set $\alpha = \dim V$.

{\bf 1.} Set $\beta = \dim V_1$. Note that $SL(\alpha)$ contains a 
subgroup naturally isomorphic
to $SL(\beta)$ and $SL(\beta)V_1+V_2$
is contained in $SL(\alpha)V$. Assuming that $V_1$
is not locally semi-simple, we get by Theorem \ref{lss-criteria}
that $\overline{SL(\beta)V_1}$ contains a representation non-isomorphic
to $V_1$, hence $\overline{SL(\alpha)V}$ contains a representation 
non-isomorphic to $V$ and $V$ is not locally semi-simple.

{\bf 2.} By Fitting's Lemma ${\rm End}(V)$ is local.
By Matsushima's criterion \cite{ma}, the stabilizer $SL(\alpha)_V$
is reductive; since ${\rm Aut}V = GL(\alpha)_V$ 
and $GL(\alpha)_V/SL(\alpha)_V$ is a subgroup in the center of 
$GL(\alpha)$,   ${\rm Aut}V$ is reductive, hence ${\rm Aut}V = {\bf k}^*$.

{\bf 3.} By {\bf 2} we know ${\rm End}(V_1) = {\bf k}$,  
${\rm End}(V_2) = {\bf k}$, and by Matsushima's criterion
${\rm End}(V)$ is reductive. The decomposition
${\rm End}(V) = \oplus_{i,j = 1,2} {\rm Hom}(V_i, V_j)$ implies
that ${\rm End}(V)$ is either ${\bf k}\oplus {\bf k}$
or ${\rm End}({\bf k}^2)$. In the latter case 
let $H_{12}$ and $H_{21}$ be the generators of
${\rm Hom}(V_1,V_2)$ and ${\rm Hom}(V_2,V_1)$, respectively.
The isomorphism ${\rm End}(V)\cong {\rm End}({\bf k}^2)$
implies $H_{21}H_{12}\in {\rm End}(V_1)$ does not vanish,
hence is a scalar operator on $V_1$.
So we have $V_1\cong V_2$, a contradiction.
\end{proof}

A representation $V$ such that ${\rm Aut}(V) = {\bf k}^*$ 
is called {\it Schurian}. The converse to \ref{nec_lss}.2 is not true, 
i.e., not any Schurian representation
is locally semi-simple:

\begin{example}
Let $Q$ be the quiver with one vertex and two attached loops.
Let $V$ be a 2-dimensional representation of $Q$ 
such that the corresponding pair of
endomorphisms of ${\bf k}^2$ generates the algebra $B$
of the upper triangular matrices, in a  basis.
Then ${\rm End}(V)$ is the centralizer of $B$ in ${\rm End}({\bf k}^2)$,
so $V$ is Schurian. Since the center of $GL(\alpha)$ acts trivially
on representations, the locally semi-simple representations
are in this case just the semi-simple representations.
A semi-simple Schurian representation must be simple,
but $V$ has a 1-dimensional subrepresentation,
so $V$ is not semi-simple.
\end{example}

Recall that each quiver $Q$ determines two forms on ${\bf Z}^{Q_0}$,
the Tits quadratic form 
$q_Q(\alpha) = \sum_{i\in Q_0} \alpha_i^2 - \sum_{\varphi\in Q_1} \alpha_{t\varphi}\alpha_{h\varphi}$,
and the Euler bilinear form: 
\begin{equation}
\langle \alpha, \beta \rangle = \sum_{i\in Q_0} \alpha_i\beta_i - \sum_{\varphi\in Q_1} \alpha_{t\varphi}\beta_{h\varphi}.
\end{equation} 
\noindent Note also that the Euler form is not symmetric and
$\langle \alpha, \alpha \rangle = q_Q(\alpha)$. 

\begin{proposition}
If $V$ is a Schurian representation and $q_Q(\dim V) = 1$ 
(in other words, $\alpha=\dim V$ is a real Schur root), then $V$ is locally semi-simple.
\end{proposition}

\begin{proof}
The hypothesis implies that the $GL(\alpha)$-orbit of $V$ is dense in 
$R(Q,\alpha)$,
so the generic stabilizer for the action of $SL(\alpha)$ is trivial. By \cite{po},
generic $SL(\alpha)$-orbit are closed. Hence, $SL(\alpha)V$ is closed.
\end{proof}


\section{Semi-invariants of quivers and perpendicular categories.}

The character group of $GL(\alpha)$ is generated by the determinants
of the $GL(\alpha_a)$-factors, $a\in Q_0$ so is isomorphic to ${\bf Z}^{Q_0}$
such that $\chi\in {\bf Z}^{Q_0}$ gives rise to the character
$\overline{\chi}=\prod_{a\in Q_0,\alpha_a>0}\det_a^{\chi_a}$.
 We also can think of 
$\chi$ as of an integer function on the dimensions of representations
such that $\chi(\alpha)=\sum_{a\in Q_0}\chi_a\alpha_a$.
We have:
\begin{equation}
{\bf k}[R(Q,\alpha)]^{SL(\alpha)} = 
\oplus_{\chi\in {\bf Z}^{Q_0}} {\bf k}[R(Q,\alpha)]^{(GL(\alpha))}_{\chi},
\end{equation}
where ${\bf k}[R(Q,\alpha)]^{(GL(\alpha))}_{\chi} = 
\{f\in {\bf k}[R(Q,\alpha)]\vert 
gf = \overline{\chi} (g)f, \forall g\in GL(\alpha)\}$.

Recall that for dimension vectors $\alpha,\beta\in {\bf Z}_+^{Q_0}$ such that
$\langle \alpha, \beta \rangle = 0$ Schofield introduced in \cite{sch} a function
$c$ on $R(Q,\alpha)\times R(Q,\beta)$ such that:

$(i)$ $c(V,W)\neq 0$ if and only if ${\rm Hom}(V,W) = 0$

$(ii)$ $c$ is $GL(\alpha)\times GL(\beta)$-semi-invariant;
if $c(V,.)\not\equiv 0$, then its character is equal $\langle \alpha,.\rangle$;
if $c(.,W)\not\equiv 0$, then its character is equal $-\langle .,\beta\rangle$.

Derksen and Weyman proved in \cite{dw} that each vector
space ${\bf k}[R(Q,\alpha)]^{(GL(\alpha))}_{\sigma}$ is generated by
the functions $c_W = c(.,W)$ such that for any $\alpha\in {\bf Z}^{Q_0}$,
$-\langle \alpha,\dim W\rangle = \sigma(\alpha)$. 
Recall the Ringel formula (\cite{ri}): 

\begin{equation}\label{ringel}
\dim {\rm Ext}(V,W) = \dim {\rm Hom}(V,W) - \langle \dim V, \dim W\rangle.
\end{equation}

\noindent This formula and the above properties imply that for a given
$V\in R(Q,\alpha)$ the semi-invariants $c_W$ such that $c_W(V)\neq 0$
correspond to the representations $W$ such that ${\rm Hom}(V,W)=0,  {\rm Ext}(V,W)=0$.
Schofield called in \cite{sch} the set of such representations, $V^{\perp}$,
the {\it right perpendicular category} of $V$. The {\it left} perpendicular
category, $^{\perp}V$, is defined similarly. 
Note that $S\in V^{\perp}$ is equivalent to $c(V,S)\neq 0$
and the same for the left category.
Schofield proved 
that the perpendicular categories are Abelian subcategories.
In particular, simple objects in $V^{\perp}$ are Schurian representations, 
homomorphisms between non-isomorphic simple objects are trivial,
any representation has a unique Jordan-H{\"o}lder decomposition.

On the other hand, we have the notion of  $\chi$-stability 
(see Remark \ref{mum}) and for representations of quivers King proved
in \cite[Proposition~3.1]{ki} that $V$ is $\chi$-stable
if and only if $\chi(\dim V) = 0$ and $\chi(\dim V') < 0$ for any 
subrepresentation $V'\subseteq V$ different from 0 and $V$. 

\begin{proposition}\label{simple_is_stable} 
$S\in ^{\perp}W$ is a simple object if and only if 
$S$ is $-\langle .,\dim W\rangle$-stable.
\end{proposition}
\begin{proof}
$S$ is a simple object means that there are no proper subrepresentations 
$S'\subseteq S$ such that ${\rm Hom}(S',W)=0$ and $\langle \dim S',\dim W\rangle = 0$.
So the "if" part follows. Now assume $S$ is simple in $^{\perp}W$
and $S'\subseteq S$ is a proper subrepresentation.
This means that $c(S,W)\neq 0$ and $c(S',W)=0$.
By \cite[Lemma~1]{dw},  $\langle\dim S',\dim W\rangle < 0$ would
contradict our condition $c(S,W)\neq 0$
and $\langle\dim S',\dim W\rangle = 0$ would imply
that $c(S',W)\neq 0$. Thus we have: $\langle\dim S',\dim W\rangle >0$.
\end{proof}

Now we give a criterion for a representation to be locally semi-simple.
The above discussion shows that it is sufficient to find out
for representations $V,W$ such that ${\rm Hom}(V,W)=0,  {\rm Ext}(V,W)=0$
whether the $GL(\dim V)$-orbit of $V$ is closed in $R(Q,\dim V)_{c_W}$.

\begin{theorem}\label{crit_lss}
The orbit of $V$ is closed in $R(Q,\dim V)_{c_W}$ if and only
if $V = p_1S_1+\cdots+p_tS_t$, where $S_1,\cdots,S_t$ are simple
objects in $^{\perp}W$.
\end{theorem}

\begin{remark}
Taking into account Proposition \ref{simple_is_stable} one can see
that this Theorem is similar to \cite[Proposition~3.2]{ki}. 
However, the property of the orbit to be closed in $R(Q,\dim V)_{c_W}$
is more strong than that to be closed in the open set of semi-stable points
as in \cite{ki}.
\end{remark}

\begin{proof}
Clearly the sum of Jordan-H{\"o}lder factors of $V$ in $^{\perp}W$
belongs to the closure of $V$-orbit and to $R(Q,\dim V)_{c_W}$.
So if the orbit is closed, then $V$ is isomorphic to the sum
of simple objects. Conversely, assume that  $V$ is a sum of simple objects
and let $U$ be the closed orbit in the closure of the orbit of $V$
in $R(Q,\dim V)_{c_W}$. Then there is a 1-parameter subgroup 
$g(t)\in GL(\dim V), t\in {\bf k}^*$ such that
$\lim_{t\to 0} g(t)V\subseteq U$. Considering the $g(t)$-eigenspace 
decomposition of $V(i), i\in Q_0$ one easily sees
that the limit exists means that these eigenspaces yield
a filtration of $V$ in $^{\perp}W$ such that $U$ is the
associated graded of $V$. Hence, the  Jordan-H{\"o}lder factors
for $U$ are the same as for $V$, so $V$ belongs to the closure
of the orbit of $U$, in other words, the orbits are equal.
\end{proof}

\begin{remark} Clearly, the Theorem implies Proposition \ref{nec_lss} above.
\end{remark}

\begin{example}\label{A_n_example}
Let $Q=A_n$: 
${\circ} \xrightarrow{} {\circ} \xrightarrow{}\cdots {\circ} \xrightarrow{} {\circ}$,
where $n$ stands for the number of vertices. Let $\varepsilon_1, \cdots, \varepsilon_n$
denote the standard basis of ${\bf Z}^n$. It is well-known that
the indecomposable representations of $A_n$ are the representations
$S_{ij}, 1 \leq i\leq j\leq n$ such that 
$\dim S_{ij} = \varepsilon_i+\cdots+\varepsilon_j$ and these representations
are Schurian. It can be directly verified ( and follows e.g. from \cite[Th.~10]{sh})
\begin{equation}\label{A_n_hom}
{\rm Hom}(S_{kl},S_{ij})\neq 0 \iff i\leq k \leq j \leq l.
\end{equation}
\noindent Hence, the condition 
${\rm Hom}(S_{ij},S_{kl})=0={\rm Hom}(S_{kl},S_{ij})$ is equivalent to:
\begin{equation}\label{A_n_hom_hom}
j < k,\; {\rm or}\; l < i,\; {\rm or}\; i<k\leq l <j,\; {\rm or} \; k<i\leq j <l,
\end{equation}
\noindent in other words the segments $[i,j]$ and $[k,l]$ 
either are disjoint sets or one of them contains another in the interior.

\begin{proposition}\label{A_n_lss}
Let $V = \oplus_{p,q}m_{pq}S_{pq}$ and set 
$I=\{(p,q)\vert m_{pq} >0\}$.
Then $V$ is locally semi-simple if and only if any pair
$(i,j),(k,l)\in I$ meets condition {\rm (\ref{A_n_hom_hom})}.
\end{proposition}

\begin{proof}
The "only if" part follows from \ref{nec_lss}. To prove the "if" part
we first observe that the condition $\langle \dim S_{kl}, \dim S_{ij}\rangle =0$
implies ${\rm Hom}(S_{kl},S_{ij})=0$, because of (\ref{A_n_hom}). 
Using condition {\rm (\ref{A_n_hom_hom})}, one can show that
the dimensions $\alpha_1,\cdots,\alpha_t$ of the representations
$S_{pq},(p,q)\in I$ are linear independent and moreover, 
the sublattice 
$$\langle \alpha_1,\cdots,\alpha_t\rangle^{\perp} = 
\{\beta\in {\bf Z}^n\vert \langle \alpha_i,\beta\rangle=0,i=1,\cdots,t\}$$

\noindent is generated by $n-t$ linear independent roots. 
By the above observation the corresponding $n-t$
indecomposable representations belong to $V^{\perp}$; 
let $W=R_1+\cdots+R_k$ be the sum of all simple factors of these. 
Clearly, $\dim R_1,\cdots,\dim R_k$ also generate
$\langle \alpha_1,\cdots,\alpha_t\rangle^{\perp}$;
since homomorphisms between simple objects $R_i,R_j$ are trivial for
$i\neq j$,  $\dim R_1,\cdots,\dim R_k$ are linearly independent,
so $k = n-t$ and $\alpha_1,\cdots,\alpha_t$
is a basis of the sublattice
$^{\perp}\langle \dim R_1,\cdots,\dim R_{n-t}\rangle$.
By construction, $V\in ^{\perp}W$. We claim that $S_{pq},(p,q)\in I$
are simple objects in $^{\perp}W$ and this implies the assertion, thanks
to Theorem \ref{crit_lss}. Indeed, assume that a summand, $S_{ij}$ 
is not simple, i.e., a proper subrepresentation $S'\subseteq S_{ij}$ 
belongs to $^{\perp}W$. Any proper subrepresentation of $S_{ij}$ is isomorphic
to $S_{kj}$ with $i<k\leq j$, so $\dim S_{kj}$ is a linear combination
of $\alpha_1,\cdots,\alpha_t$. Using condition (\ref{A_n_hom_hom}),
one can easily see this is false.
\end{proof}
\end{example}

\section{Decompositions and slices.}

Let $V$ be a locally semi-simple representation
of $Q$, $\dim V= \alpha$; by Proposition \ref{nec_lss}, we know: 
$V=\oplus_{i=1}^t m_iS_i$,
where $S_i$ are pairwise non-isomorphic Schurian representations with
trivial homomorphism spaces between them. 
Hence, ${\rm Aut}(V)\cong \prod_{i=1}^t GL(m_i)$.
Note that the group ${\rm Aut}(V)$ and its embedding to
$GL(\alpha)$ are completely determined by the decomposition
$\alpha = m_1\dim S_1 + \cdots + m_t \dim S_t$.

\begin{definition}
A decomposition $\alpha = \sum_{i=1}^t m_i\beta_i$ with 
$\alpha,\beta_1,\cdots,\beta_t\in {\bf Z}_+^{Q_0}$, $m_i\in{\bf N}$ is called
{\it locally semi-simple} 
if for each $i$ there exists a representation 
$S_i$ such that $\dim S_i = \beta_i$, $\dim {\rm Hom}(S_i,S_j) = \delta_{ij}$, 
and $V=\oplus_{i=1}^t m_iS_i$ is a locally semi-simple representation.
If moreover $V$ is generic, then we call this decomposition generic.
\end{definition}

Note that a locally semi-simple decomposition determines the isomorphism class
of the representation if and only if all the components are real Schur roots.
Note also that there can be equal summands $\beta_i=\beta_j=\beta$ 
in such a decomposition; the condition ${\rm Hom}(S_i,S_j)=0$ implies
that $\beta$ is an imaginary root. 

Assume that $V=\oplus_{i=1}^t m_iS_i$ is locally semi-simple.
By Ringel's formula (\ref{ringel}) we have:
$\delta_{ij}-\langle\dim S_i,\dim S_j\rangle=\dim {\rm Ext}(S_i,S_j) \geq 0$.
Le Bruyn and Procesi showed in \cite{lbp} that
the slice representations for $V$ semi-simple can be
expressed in terms of quivers.  
Following \cite{lbp}, we introduce a quiver
$\Sigma_V$ with vertices $a_1,\cdots,a_t$ corresponding to 
the summands $S_1,\cdots,S_t$ and $\delta_{ij} - \langle\dim S_i,\dim S_j\rangle$ arrows
from $a_i$ to $a_j$; set $\gamma = (m_1,\cdots,m_t)\in {\bf Z}^{(\Sigma_V)_0}$. 
It is known (see e.g. \cite{kr}) that for any representation
$W$ the normal space to the isomorphism class of $W$ at $W$,
$R(Q,W)/T_WGL(\dim W)W$, is isomorphic to ${\rm Ext}(W,W)$.
Hence, we get a helpful form of the slice representation 
$\sigma_V$ of $V$ (the same as in \cite{lbp} for the semi-simple case):
$(GL(\alpha)_V,\sigma_V)= ({\rm Aut}(V), {\rm Ext}(V,V)) =$ 
\begin{equation}\label{local_quiver}
(\prod_{i=1}^t GL(m_i), \bigoplus_{i,j = 1}^t {\rm Ext}(S_i,S_j)\otimes 
{\rm Hom}({\bf k}^{m_i},{\bf k}^{m_j})) = (GL(\gamma),R(\Sigma_V,\gamma)).
\end{equation}
\noindent Let $D_V:{\bf Z}^{(\Sigma_V)_0}\to {\bf Z}^{Q_0}$ denote 
the linear map taking $i$-th basis vector of ${\bf Z}^{(\Sigma_V)_0}_+$
to $\dim S_i, i=1,\cdots,t$. Note that $D_V(\gamma)=\alpha$. The definition
of the Euler form yields:
\begin{proposition}\label{isometry}
$D_V$ preserves the Euler form:
$\langle D_V(\gamma_1),D_V(\gamma_2)\rangle = \langle\gamma_1,\gamma_2\rangle$.
\end{proposition}

\begin{proposition}\label{slice_decomp}
Consider a decomposition $\gamma = \sum_{j=1}^s p_j\rho_j$ and
the corresponding decomposition $\alpha = \sum_{j=1}^s p_j D_V(\rho_j)$
of $\alpha$.

{\bf 1.} If the decomposition of $\gamma$ is generic, then that of $\alpha$ is.

{\bf 2.} If the decomposition of $\gamma$ is locally semi-simple, 
then that of $\alpha$ is.

{\bf 3.} If the decomposition of $\gamma$ is generic locally semi-simple, 
then that of $\alpha$ is.
\end{proposition}
\begin{proof}
A general remark is that by Luna's slice theorem for
any representation $W\in R(\Sigma_V,\gamma)$ there exists
a representation $V_1\in R(Q,\alpha)$ with ${\rm Aut}(V_1)={\rm Aut}(W)$,
where ${\rm Aut}(W)\subseteq GL(\gamma)$ is embedded to $GL(\alpha)$ via the embedding
$GL(\gamma)={\rm Aut}(V)\subseteq GL(\alpha)$. 
Therefore if the maximal torus of ${\rm Aut}(W)$
corresponds to the given decomposition of $\gamma$, then
the maximal torus of ${\rm Aut}(V_1)$
corresponds to the given decomposition of $\alpha$.
Now {\bf 1} follows from the fact that the generic decompositions are determined
by generic stabilizers and by Luna's slice theorem if $W$ is generic for
$(GL(\alpha)_V,\sigma_V)$, then $V_1$ is generic for $(GL(\alpha),R(Q,\alpha))$.
By Proposition \ref{lss_from_slice} $W$ is locally semi-simple
implies $V_1$ is locally semi-simple, so we proved {\bf 2}.
Applying Corollary \ref{generic_lss_from_slice}, we also get {\bf 3}.
\end{proof}

\begin{proposition}\label{D_V_prop}
{\bf 1.}
If $\gamma_1\in {\bf Z}^{(\Sigma_V)_0}_+$ is a (real) Schur root, 
then $D_V(\gamma_1)$ is.  

{\bf 2.} If $\gamma_1,\gamma_2\in {\bf Z}^{(\Sigma_V)_0}_+$ are real Schur roots,
$W_1\in R(\Sigma_V,\gamma_1),W_2\in R(\Sigma_V,\gamma_2)$
and $V_1\in R(Q,D_V(\gamma_1)),V_2\in R(Q,D_V(\gamma_2))$ are Schurian 
representations, then 
${\rm Ext}(W_1,W_2)=0={\rm Ext}(W_2,W_1)$, 
implies ${\rm Ext}(V_1,V_2)=0={\rm Ext}(V_2,V_1)$ and 
${\rm Hom}(W_1,W_2)=0={\rm Hom}(W_2,W_1)$, 
implies ${\rm Hom}(V_1,V_2)=0={\rm Hom}(V_2,V_1)$.
\end{proposition}
\begin{proof}
Note that the map $D_V$ depends on the indecomposable summands of $V$
not of $V$ itself.  So in {\bf 1} we may assume that $\dim V= D_V(\gamma_1)$.
Then by \ref{slice_decomp}.1 $\dim V = \dim V$ is the generic decomposition,
that is $D_V(\gamma_1)$ is a Schur root. By \ref{isometry},
$q_Q(D_V(\gamma_1)) =q_{\Sigma_V}(\gamma_1)$. So $\gamma_1$
is real implies $D_V(\gamma_1)$ is. In {\bf 2} we assume
$\dim V =\alpha= D_V(\gamma_1)+ D_V(\gamma_2)$. Then by \cite{kac} the condition
${\rm Ext}(W_1,W_2)=0={\rm Ext}(W_2,W_1)$ implies
that $\gamma =\gamma_1+\gamma_2$ is the generic decomposition.
Then by \ref{slice_decomp}.1, $\alpha =  D_V(\gamma_1)+ D_V(\gamma_2)$
is the generic decomposition so again applying \cite{kac} we get
${\rm Ext}(V_1,V_2)=0={\rm Ext}(V_2,V_1)$.
The condition ${\rm Hom}(W_1,W_2)=0={\rm Hom}(W_2,W_1)$
yields ${\rm Aut}(W_1+W_2)$ is the corresponding embedding
of $({\bf k}^*)^2$ to $GL(\gamma)$. By Luna's slice theorem
there exists a representation $V'\in R(Q,\alpha)$ with 
${\rm Aut}(V')$ being the image of ${\rm Aut}(W_1+W_2)$
under the embedding $GL(\gamma)\subseteq GL(\alpha)$. 
This means that $V'=V'_1+V'_2$ with $V'_1,V'_2$
indecomposable of dimensions $D_V(\gamma_1),D_V(\gamma_2)$
and ${\rm Hom}(V'_1,V'_2)=0={\rm Hom}(V'_2,V'_1)$.
Clearly, this is equivalent to what we assert.
\end{proof}

Now apply the Luna-Richardson theorem to this situation.
Assume that $\alpha = \sum_{i=1}^t m_i\beta_i$ is a generic
locally semi-simple decomposition. For a locally semi-simple
decomposition we observed above that equal summands $\beta_i=\beta_j$ 
must be imaginary roots. Besides in generic case the multiplicity $m_j$ of 
any imaginary root $\beta_j$ is equal 1. Indeed, if we have
a locally semi-simple representation $V=\oplus_{i=1}^t m_iS_i$
then by Theorem \ref{crit_lss} $S_1,\cdots,S_t$ are simple
objects in the category $^{\perp}W$ for some $W$.
Clearly, being perpendicular to $W$ and being simple object
in $^{\perp}W$ are open conditions on $R(Q,\beta_j)$. 
Since $R(Q,\beta_j)$ contains infinitely many isomorphism classes
of indecomposable representations, we could replace 
$m_jS_j$ with $m_j>1$ by a sum of $m_j$ generic representations of dimension
$\beta_j$ to get a locally semi-simple representation
with a smaller automorphism group.

\begin{theorem}\label{lr_quiver}
{\bf 1.} A generic locally semi-simple decomposition has the form:
\begin{equation}\label{glssd}
\alpha = \underbrace{\delta_1+\cdots+\delta_1}_{p_1\; summands}+\cdots+
\underbrace{\delta_r+\cdots+\delta_r}_{p_r\; summands}+
m_1\beta_1+\cdots+m_s\beta_s,
\end{equation}
where $\delta_1,\cdots,\delta_r$ are pairwise non-equal
imaginary Shur roots and $\beta_1,\cdots,\beta_s$ are pairwise non-equal
real Shur roots.

{\bf 2.} Generic stabilizer $H$ of a locally semi-simple point 
is isomorphic to $\prod_{i=1}^r({\bf k}^*)^{p_i}$
$\times \prod_{j=1}^s GL(m_j)$.
The linear group $(N_{GL(\alpha)}(H)/H,R(Q,\alpha)^{H})$ is isomorphic to
\begin{equation}\label{weyl_group}
\bigoplus_{i=1}^r \underbrace{((GL(\delta_i),R(Q,\delta_i))\oplus\cdots
\oplus((GL(\delta_i),R(Q,\delta_i))}_{S_{p_i}\;permutes\;p_i\;summands}
\oplus\bigoplus_{j=1}^s (GL(\beta_j), R(Q,\beta_j)).
\end{equation}

{\bf 3.} ${\bf k}[R(Q,\alpha)]^{SL(\alpha)}\cong$
${\bf k}[R(Q,p_1\delta_1)\oplus\cdots R(Q,p_r\delta_r)\oplus
         R(Q,\beta_1)\oplus\cdots R(Q,\beta_s)]^G$,
where $G\subseteq GL(p_1\delta_1)\times\cdots\times GL(p_r\delta_r)\times 
                  GL(\beta_1)\times\cdots\times GL(\beta_s)$ 
consists of the elements such that for each vertex the product of determinants
is 1. Moreover, generic
$G$-orbit in 
$R(Q,p_1\delta_1)\oplus\cdots R(Q,p_r\delta_r)\oplus
 R(Q,\beta_1)\oplus\cdots R(Q,\beta_s)$ is closed.

{\bf 4.} ${\bf k}[R(Q,\alpha)]^{SL(\alpha)}\cong$
$\oplus_{\chi\in\Lambda}
{\bf k}[R(Q,p_1\delta_1)]^{(GL(p_1\delta_1))}_{\chi}
\otimes\cdots\otimes
{\bf k}[R(Q,p_r\delta_r)]^{(GL(p_r\delta_r))}_{\chi}$,
where $\Lambda=\{\chi\in{\bf Z}^{Q_0}\vert 
{\bf k}[R(Q,\beta_j)]^{(GL(\beta_j))}_{\chi}\neq 0,
j=1,\cdots,s\}$.
\end{theorem}

\begin{proof}
The assertion {\bf 1} is showed above. 
The form of $H$ in {\bf 2} follows from (\ref{glssd}).
For any vertex $a\in Q_0$ each summand $\rho=\delta_i$ or $\beta_j$
yields an isotypical component of the $H$-module ${\bf k}^{\alpha_a}$
being the sum of $\rho_a$ irreducible factors of type
$(GL_m,{\bf k}^m)$, where $m=1$ for $\delta_i$ and $m=m_j$ for $\beta_j$.
These isotypical components are stable with respect to the centralizer 
$Z_{GL(\alpha)}(H)$ of $H$ and each of them yields
a factor $(GL(\rho),R(Q,\rho))$ of $(Z_{GL(\alpha)}(H)/H,R(Q,\alpha)^H)$.
Elements of $N_{GL(\alpha)}(H)\setminus Z_{GL(\alpha)}(H)$
induce an outer automorphism of the group $H$ and the corresponding
permutation of the isotypical components in 
each space ${\bf k}^{\alpha_a}$. 
Clearly, the isotypical 
components corresponding to non-equal summands can not be
permuted, so $N_{GL(\alpha)}(H)$ is contained
in the extension of $Z_{GL(\alpha)}(H)$ by the groups $S_{p_i}, i=1,\cdots,r$.
and the latter extension does normalize $H$ so {\bf 2} is proved.

By Proposition \ref{lr_lss}, 
${\bf k}[R(Q,\alpha)]^{SL(\alpha)}\cong 
{\bf k}[R(Q,\alpha)^H]^{N_{SL(\alpha)}(H)}$. 
Consider a subgroup $N\triangleleft N_{GL(\alpha)}(H)$ consisting
of the elements such that the restrictions  
to the irreducible $N_{GL(\alpha)}(H)$-submodule of ${\bf k}^{\alpha_a}$
are unimodular for any $a\in Q_0$. 
By (\ref{weyl_group}) $N$ acts independently
on the summands of $R(Q,\alpha)^H$ and we have:
${\bf k}[R(Q,\alpha)^H]^N\cong$  
$(\otimes_{i=1}^r{\bf k}[R(Q,\delta_i)\oplus\cdots\oplus R(Q,\delta_i)]^N)
\otimes (\otimes_{j=1}^s {\bf k}[R(Q,\beta_j)]^N)$.
Note that $p_i\delta_i=\delta_i+\cdots+\delta_i$ is the generic locally
semi-simple decomposition of $p_i\delta_i$.
Note also that $N_{GL(\alpha)}(H)$ acts on 
$R(Q,\delta_i)\oplus\cdots\oplus R(Q,\delta_i)$
as $N$ extended by the center of $GL(p_i\delta_i)$.
Therefore by \ref{lr_lss},
${\bf k}[R(Q,\delta_i)\oplus\cdots\oplus R(Q,\delta_i)]^N$ 
and ${\bf k}[R(Q,p_i\delta_i)]^{SL(p_i\delta_i)}$
are isomorphic as $\Xi(GL(p_i\delta_i))$-graded algebras.

Next, fix $j\in\{1,\cdots,s\}$, set $m=m_j,\beta=\beta_j$ and consider 
the direct summand of $R(Q,\alpha)^H$ corresponding to $m\beta$.
The determinant $\det_a$ restricts to $GL(\beta)\subseteq N_{GL(\alpha)}(H)$ 
as the $m$-th power of the corresponding determinant
on $GL(\beta_a)$. Therefore $N$ acts on $R(Q,\beta)$
as the group $SL(\beta)\Gamma_m$, where 
$\Gamma_m=\{g\in GL(\beta)\vert g_a = c_a{\rm Id}, c_a\in \sqrt[m]{1}\}$
is a finite group. Since $\beta$ is a real Schur root,
${\bf k}[R(Q,\beta)]^{SL(\beta)}$ is generated by semi-invariants
with linear independent weights such that their common kernel
is the group ${\bf k}^*$ of scalar operators. One can deduce from this that
for each weight $\chi$ there is an element $w\in\Gamma_m$ such that 
$\chi(w)$ is a prime unity root of order $m$ and all other
weights take $w$ to 1. Consequently, ${\bf k}[R(Q,\beta)]^N$
is generated by the $m$-th powers of the generators of 
${\bf k}[R(Q,\beta)]^{SL(\beta)}$. So we have an isomorphism
of ${\bf k}[R(Q,\beta)]^{SL(\beta)}$ onto ${\bf k}[R(Q,\beta)]^N$
taking $GL(\beta)$-eigenvectors of weight $\chi$ to 
$N_{GL(\alpha)}(H)/N$-eigenvectors of weight $m\chi$.

Thus we proved an isomorphism
\begin{equation}
{\bf k}[R(Q,\alpha)^H]^N\cong\bigotimes_{i=1}^r
{\bf k}[R(Q,p_i\delta_i)]^{SL(p_i\delta_i)}\otimes
\bigotimes_{j=1}^s{\bf k}[R(Q,\beta_j)]^{SL(\beta_j)},
\end{equation}
where a subspace of weight 
$\sigma_1+\cdots+\sigma_r+\chi_1+\cdots+\chi_s$ in the
right hand side algebra corresponds to that of weight
$\sigma_1+\cdots+\sigma_r+m_1\chi_1+\cdots+m_s\chi_s$ in 
${\bf k}[R(Q,\alpha)^H]^N$. Consequently, the 
$GL(p_1\delta_1)\times\cdots\times GL(\beta_s)$-weights vanishing
on $G$ correspond to the $N_{GL(\alpha)}(H)/N$-weights
vanishing on $N_{SL(\alpha)}(H)$ so the subalgebra of
$N_{SL(\alpha)}(H)$-invariants corresponds under the isomorphism
to that of $G$-invariants. By \ref{lr_lss} generic
$N_{SL(\alpha)}(H)$-orbits are closed in $R(Q,\alpha)^H$;
clearly, this implies the same for the $G$-orbits,
 and {\bf 3} is proved.

Elements of
$(\otimes_{i=1}^r {\bf k}[R(Q,p_i\delta_i)]^{(GL(p_i\delta_i))}_{\sigma_i})
\otimes 
(\otimes_{j=1}^s {\bf k}[R(Q,\beta_j)]^{(GL(\beta_j))}_{\chi_j})$
are $G$-inva\-riant if and only
if the corresponding character
$\sigma_1+\cdots+ \sigma_r + \chi_1 +\cdots + \chi_s$
is a linear combination of the sums of determinants for each vertex,
that is, $\sigma_1 =\sigma_2 =\cdots =\sigma_r = \chi_1 =\cdots = \chi_s$
as elements of ${\bf Z}^{Q_0}$.
So the algebra of $G$-invariants is equal to
$\oplus_{\chi\in\Lambda}
(\otimes_{i=1}^r {\bf k}[R(Q,p_i\delta_i)]^{(GL(p_i\delta_i))}_{\chi})
\otimes
(\otimes_{j=1}^r {\bf k}[R(Q,\beta_j)]^{(GL(\beta_j))}_{\chi})$.
On the other hand, $\beta_j$ are real Schur roots, hence,
$\dim {\bf k}[R(Q,\beta_j)]^{(GL(\beta_j))}_{\chi}=1$ for any 
$j$ and $\chi\in\Lambda$. So restricting $G$-invariants
to $\oplus_{i=1}^r R(Q,p_i\delta_i)$, we get an 
isomorphism of ${\bf k}[(\oplus_{i=1}^r R(Q,p_i\delta_i))\oplus
(\oplus_{j=1}^s R(Q,\beta_j))]^G$ onto
$\oplus_{\chi\in\Lambda}
\otimes_{i=1}^r {\bf k}[R(Q,p_i\delta_i)]^{(GL(p_i\delta_i))}_{\chi}$.
\end{proof}

\section{Decompositions for $A_n$ quiver.}
In this section we describe generic and generic locally semi-simple
decompositions for $Q$ being the equioriented $A_n$-quiver that we considered
in Example \ref{A_n_example}. Since $Q$ is finite, there is a dense
isomorphism class in $R(Q,\alpha)$ for all $\alpha$ so $V$ is generic is equivalent
to $V$ having the dense orbit or ${\rm Ext}(V,V)=0$. So we are looking
for a sum of $S_{ij}$ with trivial ${\rm Ext}$-spaces for summands.
Using Ringel formula (\ref{ringel}) and (\ref{A_n_hom}),
we see that the condition 
${\rm Ext}(S_{ij},S_{kl})=0={\rm Ext}(S_{kl},S_{ij})$ is equivalent to:
\begin{equation}\label{A_n_ext_ext}
j < k-1,\; {\rm or}\; l < i-1,\; {\rm or}\; i\leq k\leq l \leq j,\; 
{\rm or} \; k\leq i\leq j \leq l.
\end{equation}
so either the distance between the segments $[i,j]$ and $[k,l]$ is
at least 2, or one of them contains another.
This property yields a simple algorithm for calculating generic decomposition,
exactly the same as in \cite[Lemma~3.3]{ri}:

\begin{algorithm}\label{alg_gen_A_n}
For $\alpha \in {\bf Z}^{Q_0}_+$ set $m=min\{\alpha_a\vert a\in Q_0\}$.
If $m>0$, then $\alpha = m(1,\cdots,1)+\pi,\pi \in {\bf Z}^{Q_0}_+$ 
and the generic decomposition of $\alpha$ is 
$\alpha = m(1,\cdots,1)+$ the terms of  the generic decomposition 
of $\pi$. 
Otherwise, if $\alpha_t=0$, then 
$\alpha = \pi+\sigma$, $\pi=(\alpha_1,\cdots, \alpha_{t-1},0,\cdots,0)$,
$\sigma = (0,\cdots,0,\alpha_{t+1},\cdots, \alpha_n)$ and
the generic decomposition of $\alpha$ is that of $\pi+$ that of $\sigma$
for the appropriate proper subquivers.
\end{algorithm}

Now we consider locally semi-simple decompositions. The following observation
follows from Proposition \ref{A_n_lss}:

\begin{proposition}\label{min_split}
If $0<m=\alpha_t<\alpha_i$ for any $i\neq t$, then any locally semi-simple
representation $V$ of dimension $\alpha$ decomposes as $V= mS_{tt}+$ other summands.
\end{proposition}

\begin{algorithm}\label{alg_gen_lss_A_n}
For $\alpha \in {\bf Z}^{Q_0}_+$ set $m=min\{\alpha_a\vert a\in Q_0\}$,
$t=min\{a\in Q_0\vert \alpha_a=m\}$, $s=max\{a\in Q_0\vert \alpha_a=m\}$.
If $m>0$, then $\alpha = \pi+\sigma+\rho+ m\dim S_{ts}$,
where $\pi = (\alpha_1,\cdots, \alpha_{t-1},0,\cdots,0)$,
$\sigma = (0,\cdots,0,\alpha_{s+1},\cdots, \alpha_n)$,
$\rho\in \langle\varepsilon_{t+1},\cdots, \varepsilon_{s-1}\rangle$
and the generic locally semi-simple decomposition of $\alpha$
is $\alpha = m\dim S_{ts}+$ the terms of the decompositions
of $\pi,\sigma,\rho$ for the appropriate proper subquivers.
Otherwise, if $\alpha_t=0$, then 
$\alpha = \pi+\sigma$, $\pi=(\alpha_1,\cdots, \alpha_{t-1},0,\cdots,0)$,
$\sigma = (0,\cdots,0,\alpha_{t+1},\cdots, \alpha_n)$ and
the generic locally semi-simple decomposition of $\alpha$ is 
that of $\pi+$ that of $\sigma$
for the appropriate proper subquivers.
\end{algorithm}

\begin{proof}
The second case $m=0$ is obvious. In the first case
set $\mu=\alpha -  m\dim S_{ts}$ and
take a representation 
$V= mS_{ts} + \mu_1 S_{11}+\cdots+ \mu_n S_{nn}, \dim V = \alpha$.
Since $\mu_t = \mu_s = 0$, by \ref{A_n_lss} $V$ is locally semi-simple.
The ${\rm Ext}$-spaces for the summands of $V$ are non-zero and 
one-dimensional only for 
${\rm Ext}(S_{ii},S_{i+1i+1})$, ${\rm Ext}(S_{t-1t-1},S_{ts})$,
and ${\rm Ext}(S_{ts},S_{s+1s+1})$.
So the graph $\Sigma_V$ is the disjoint union of $A_{n-s+t}$ on the vertices
corresponding to $S_{11},\cdots,S_{t-1t-1},S_{ts},S_{s+1s+1},S_{nn}$ and
$A_{s-t-1}$ (if $s-t\geq 2$) on the vertices
corresponding to $S_{t+1t+1},\cdots,S_{s-1s-1}$. The induced dimension $\gamma$
is $(\pi,m,\sigma)$ on $A_{n-s+t}$ and $\rho$ on $A_{s-t-1}$.
By Proposition \ref{slice_decomp}.3, the generic locally semi-simple
decomposition for $\alpha$ is the sum of that for $\rho$ and that
for $(\pi,m,\sigma)$. By Proposition \ref{min_split},
the generic  locally semi-simple decomposition for $(\pi,m,\sigma)$
is $m\varepsilon_t+$ the sum of the decompositions for $\pi$ and $\sigma$.
Applying the map $D_V$ to the summands of the decomposition, we conclude the proof. 
\end{proof}
 

\section{ Regular representations of tame quivers.}\label{section_regular}

The tame quivers can be described by several equivalent
conditions; in particular, these are the quivers with the underlying
graph being an extended Dynkin diagram of type 
$\widetilde{A_n},\widetilde{D_n},\widetilde{E_6},\widetilde{E_7},
\widetilde{E_8}$  (the number of vertices is in all cases the subscript + 1).
So let $Q$ be a tame quiver and assume additionally that $Q$
does not have oriented cycles (this is a restriction only for
the underlying graph being $\widetilde{A_n}$).

For quivers without oriented cycles Bernstein, Gelfand, and Ponomarev introduced
in \cite{bgp} {\it Coxeter functors} $C^+$ and $C^-$ 
(defined not uniquely) acting on representations
of $Q$. The corresponding linear Coxeter transformation $c$ 
is defined by 
the rule $c(\dim V)=\dim C^+ V$ for a representation
$V$ of dimension $\alpha$; note that $\dim C^-(V) = c^{-1}\dim V$.
Indecomposable representations $V$ such that $C^{+n}V=0$
for natural $n$ are called {\it preprojective}, the
{\it preinjective} representations being defined symmetrically.
Representation having neither preprojective nor preinjective
direct summands are called {\it regular}.

For tame quivers regular indecomposable representations $V$
can be described in term of a certain {\it defect} function $\sigma$
such that $V$ is regular if and only if $\sigma(\dim V)=0$.
This $\sigma$ is presented explicitly in \cite{ri} (for special
orientations), and one can easily check in all cases: 
\begin{equation}\label{sigma_delta}
\sigma(\alpha) = \langle \alpha,\delta\rangle,
\end{equation}
where $\delta$ is the non-divisible imaginary root such that 
$\langle \delta,\delta\rangle=0$.

In \cite{ri} Ringel proved that the regular representations form
an Abelian subcategory ${\mathcal R}$ closed under direct sums, direct summands,
homomorphisms, extensions etc. Note that by definition and (\ref{sigma_delta}), 
the simple regular objects are the $\delta$-stable representations. 
These simple objects are as follows. In dimension $\delta$ there is a 1-parameter
family of simple regular objects; these representations 
are called {\it homogeneous}. 
We follow \cite{ri} and denote
by $I$ the set of the dimensions of regular simple objects
different from $\delta$ and by $e_i$ the dimension corresponding to $i\in I$. 
It is known that the set $I$ consists of real Schur roots so
there is a unique simple representation $E_i$ of dimension $e_i$,
up to isomorphism.
Furthermore, the set $I$ is finite and stable with respect to the 
Coxeter transformation $c$; moreover, $c$ has at most 3 orbits 
in $I$. The sum of dimensions over a $c$-orbit is equal $\delta$.

The category ${\mathcal R}$ is connected with perpendicular ones:

\begin{proposition}\label{regular_is_perp}
Let $S\in {\mathcal R}$ be a homogeneous simple object
and let $V\in {\mathcal R}$ be an indecomposable representation.
If not all Jordan-H{\"o}lder factors of $V$ are isomorphic to $S$,
then $V\in ^{\perp}S$ and $V\in S^{\perp}$.
\end{proposition} 

\begin{proof}
By definition and (\ref{sigma_delta}), 
$\langle\dim V,\dim S\rangle = \langle\dim V,\delta\rangle = \sigma(\dim V)=0$.
On the other hand, for any $\alpha$, 
$\langle\alpha,\delta\rangle + \langle\delta,\alpha\rangle
=q_Q(\alpha+\delta) - q_Q(\alpha) - q_Q(\delta) = 0$,
because $\delta$ is in the kernel of $q_Q$.
So we need to check: ${\rm Hom}(V,S) = 0 = {\rm Hom}(S,V)$ 
or (by the Ringel formula)
${\rm Ext}(V,S) = 0 = {\rm Ext}(S,V)$. Assume the converse 
and apply induction on $\dim V$.
If ${\rm Hom}(V,S)\neq 0$, then we get 
an exact sequence $0\to V'\to V\to S\to 0$, because $S$ is simple.
So $V'$ also has a Jordan-H{\"o}lder factor 
different from $S$, since $V$ has. Consequently,
at least one of the direct summands
$V_1'$ of $V'$ meets this condition and ${\rm Ext}(S,V_1')=0$ by induction.
Decompose $V'$ as $V'=V_1'+V_2'$.
Applying the definition of the ${\rm Ext}$-functor,
one can find a subrepresentation $V_3'\subseteq V$ containing
$V_2'$ such that $V=V_1'+V_3'$. This is a contradiction,
because $V$ is indecomposable.
Analogously, ${\rm Hom}(S,V)\neq 0$ and the induction imply 
${\rm Ext}(V_1',S)=0$ for an indecomposable summand $V_1'\subseteq V/S$
and this also contradicts to $V$ being indecomposable.
\end{proof} 

Denote by ${\mathcal D}_r$ the dimensions of regular representations.
If $\alpha\notin {\mathcal D}_r$, then by \cite[Theorem~3.2]{ri}
$R(Q,\alpha)$ contains a dense orbit. Otherwise, if $\alpha\in {\mathcal D}_r$,
then $\alpha$ decomposes as $\alpha = p\delta+\sum_{i\in I}p_ie_i$
and there is a unique decomposition of such a type
with an additional condition that for every $c$-orbit there is
an element $j$ such that $p_j=0$. Ringel called this decomposition
{\it canonical}. This decomposition yields locally semi-simple
representations:

\begin{proposition}\label{canon_yield_lss}
Let $\alpha = p\delta+\sum_{i\in I}p_ie_i$
be the canonical decomposition of $\alpha\in {\mathcal D}_r$. 
Consider a representation 
$V = S_1+\cdots+S_p+\sum_{i\in I}p_iE_i$, where $S_1,\cdots,S_p$
are homogeneous representations.
Then $V$ is a locally semi-simple representation.
\end{proposition}

\begin{proof}
Take $S$ to be a homogeneous simple object non-isomorphic
to $S_1,\cdots,S_p$. By Proposition \ref{regular_is_perp}
$S_1,\cdots,S_p$ and $E_i$ for all $i\in I$ belong to $^{\perp}S$.
Since all these are $\delta$-stable, these are also
simple in $^{\perp}S$, by Proposition \ref{simple_is_stable}.
So the assertion follows from Theorem \ref{crit_lss}. 
\end{proof}

In (\ref{local_quiver}) we described the slice at a locally semi-simple
point $V$ in terms of the quiver $\Sigma_V$ with dimension $\gamma$. 
For $V$ being as in \ref{canon_yield_lss} $\Sigma_V$ has a simple structure.
Denote by $E(Q)$ the quiver with $E(Q)_0=I$ and an arrow
from $i$ to $j$ for each pair $(i,j)$ such that $c(e_i)=e_j$.
Note that $E(Q)$ is a disjoint union of circular quivers. 

\begin{proposition}\label{local_quiver_at_canon}
Let $\alpha$ and $V$ be as in Proposition \ref{canon_yield_lss}
such that $S_1,\cdots,S_p$ are pairwise non-isomorphic.
Then $(\Sigma_V,\gamma)$ is a disjoint union of 1-dimensional representations
of the loops sitting at the vertices corresponding to $S_1,\cdots,S_p$
and $(E(Q),p_i,i\in I)$.
\end{proposition}
\begin{proof}
The quiver $\Sigma_V$ is defined in terms of the Euler form 
or the ${\rm Ext}$-spaces for the summands of $V$. 
Since $\langle \delta, \delta\rangle = 0$ and
$\langle \delta, e_i\rangle =\langle e_i,\delta\rangle= 0$,
each of the vertices corresponding to $S_1,\cdots,S_p$
is incident to the unique arrow-loop and the dimension sitting
there is 1. 
It remains to describe ${\rm Ext}(E_i,E_j)$.
Applying the formula: 
\begin{equation}
\dim {\rm Ext}(U,W) = \dim {\rm Hom} (W,C^+U).
\end{equation}
(see e.g. \cite[p.219]{ri}) we get: $\dim {\rm Ext}(E_i,E_j) = \dim {\rm Hom}(E_j,C^+E_i)$.
Since $\dim{\rm Hom}(E_i,E_j)=\delta_{ij}$, $\dim {\rm Ext}(E_i,E_j)$
is either 0 or 1, the latter being equivalent to
$c(e_i)=e_j$.
\end{proof}

Now we have 3 ingredients that allow to calculate the generic and the generic 
locally semi-simple decompositions for $\alpha\in {\mathcal D}_r$.
First, given the canonical decomposition of $\alpha$, we have
a locally semi-simple representation $V$ and the description of $\Sigma_V$
in Proposition \ref{local_quiver_at_canon}. Thanks to the condition that
$p_i=0$ for at least one $i$ in each $c$-orbit,
the group $(GL(\gamma), R(\Sigma_V,\gamma))$ is isomorphic,
up to a $p$-dimensional invariant subspace to a
direct sum of groups  $(GL(\gamma_i), R(A_{n_i},\gamma_i))$.
Secondly, Proposition \ref{slice_decomp} reduces both
decompositions to the same for the quivers $A_n$.
Thirdly, Algorithms \ref{alg_gen_A_n} and \ref{alg_gen_lss_A_n} 
yield both decompositions for $A_n$.

In what concerns the generic decomposition our algorithm 
recovers that by Ringel from \cite[Theorem~3.5]{ri}.
It should be noted, however, that Ringel used an
equivalence of categories instead of the slice theorem.  

\begin{example}
Consider the quiver $Q=\widetilde{E_6}$ (over each vertex we placed the index):
$$
\begin{array}{ccccccccc}
\widetilde{E_6} :
\overset{1}{\circ} & \longrightarrow & \overset{2}{\circ} & \longrightarrow &
\overset{7}{\circ} & \longleftarrow & \overset{4}{\circ} & \longleftarrow & \overset{3}{\circ} \\
&&&&\uparrow &&&&\\
&& \overset{5}{\circ} & \longrightarrow &
\overset{6}{\circ} &&&&
\end{array}
$$

We have $\delta = (1,2,1,2,1,2,3)$ so that
$\sigma(\alpha)=3\alpha_7-\alpha_1-\cdots-\alpha_6$. 
The sequence of the vertices
in the order defined by the indices is {\it admissible}
in the sense of \cite{bgp}, i.e., for any arrow $\varphi$ holds
$h\varphi >t\varphi$. So the composition $C^+=R^+_1R^+_2\cdots R^+_7$
of the reflection functors at {\it sinks} is well-defined. 
Hence we have $c=r_1r_2\cdots r_7$ where $r_i$ is the 
reflection at the vertex $i$.
There are 3 $c$-orbits of dimensions of simple regular representations:
$e_3\to e_2\to e_1\to e_3$, $e_6\to e_5\to e_4\to e_6$, 
$e_8\to e_7\to e_8$:
\begin{equation}
e_1=(1,1,0,1,0,0,1);e_2=(0,0,1,1,0,1,1);e_3=(0,1,0,0,1,1,1);
\end{equation}

\begin{equation}
e_4=(1,1,0,0,0,1,1);e_5=(0,0,0,1,1,1,1);e_6=(0,1,1,1,0,0,1);
\end{equation}

\begin{equation}
e_7=(0,1,0,1,0,1,1);e_8=(1,1,1,1,1,1,2).
\end{equation}

For example take $\alpha = (6,10,7,14,5,9,17)$. The canonical
decomposition of $\alpha$ is: $\alpha = 2\delta+ 3e_1+2e_2+
2e_5+ 2e_6+ e_8$. So $(E(Q),\gamma)$ is the direct sum
of $(A_2,(2,3))$, $(A_2,(2,2))$, and $(A_1,(1))$.
Applying Algorithm \ref{alg_gen_A_n}, we get the generic decomposition
for $(E(Q),\gamma)$:  $(A_2,2(1,1)+(0,1))$, $(A_2,2(1,1))$, and  $(A_1,(1))$.
So by Proposition \ref{slice_decomp}.1, the generic decomposition
of $\alpha$ is $\alpha = 2\delta+ 2(e_1+e_2)+ e_1 + 2(e_5+ e_6) + e_8$,
where $e_1+e_2$ and $e_5+e_6$ are real Schur roots.
Next, applying Algorithm \ref{alg_gen_lss_A_n}, we get the generic 
locally semi-simple decomposition
for $(E(Q),\gamma)$:  $(A_2,2(1,0)+3(0,1))$, $(A_2,2(1,1))$, and  $(A_1,(1))$.
So by Proposition \ref{slice_decomp}.3, the generic locally semi-simple
decomposition
of $\alpha$ is $\alpha = 2\delta+ 3e_1 + 2e_2 + 2(e_5+ e_6) + e_8$. 
\end{example}

\section{Semi-invariants of tame quivers.}

The algebras of semi-invariants 
of tame quivers $Q$ have been studied in several papers
including \cite{ri},\cite{hh},\cite{schw}. In \cite{skw} Skowronsky
and Weyman proved that ${\bf k}[R(Q,\alpha)]^{SL(\alpha)}$ is a 
{\it complete intersection} for any $\alpha$; 
moreover in most cases ${\bf k}[R(Q,\alpha)]^{SL(\alpha)}$ 
is a polynomial algebra and in all other cases is a hypersurface.

Note that after \cite{kac} it is known that the reflection functors 
give rise to so called {\it castling transforms} of semi-invariants,
so given a description of semi-invariants for $Q$ and $\alpha$,
one can describe the semi-invariants for any quiver and dimension
obtained by reflection functors. In particular, one may fix
a convenient orientation for $Q$ (in the case of $\widetilde{A_n}$,
one of the convenient orientations). If $\alpha\notin{\mathcal D}_r$,
then by \cite[Theorem~3.2]{ri}, $R(Q,\alpha)$ contains a dense orbit,
hence ${\bf k}[R(Q,\alpha)]^{SL(\alpha)}$ is a polynomial algebra
by the theorem of Sato-Kimura (\cite{sk}). Moreover, one can always apply one of 
the Coxeter functors $C^+$ or $C^-$ and describe the semi-invariants
of $Q$ in dimension $\alpha$ in terms of the castling transforms
of those in dimension $\beta=c(\alpha)$ or $c^{-1}(\alpha)$, respectively.
It is well known that for $\alpha\notin{\mathcal D}_r$ this process is
not cyclic and in the end we reduce the question to  
$\alpha$ being the dimension of a representation of a projective or an injective
module where the semi-invariants are obvious (see an example of such an
approach in \cite{schw} for $\widetilde{D_4}$ quiver).
That is why we may and will assume from now on: 
$\alpha = p\delta+\sum_{i\in I}p_ie_i\in{\mathcal D}_r$.

Ringel described the field ${\bf k}(R(Q,\alpha))^{GL(\alpha)}$ 
of invariants. Namely, he constructed semi-invariants 
$f_0,\cdots,f_p$ of weight $\sigma$ and proved in 
\cite[Theorem~4.1]{ri} that the fractions 
$\frac{f_1}{f_0},\cdots,\frac{f_p}{f_0}$
generate ${\bf k}(R(Q,p\delta))^{GL(p\delta)}$.
Moreover, it is stated on \cite[p.237]{ri} that 
$f_0,\cdots,f_p$ form a basis of 
${\bf k}[R(Q,p\delta)]^{(GL(p\delta))}_{\sigma}$
and one can actually deduce this from the proof of \cite[Theorem~4.1]{ri}.

First consider the homogeneous case $\alpha= p \delta$.
The generators of ${\bf k}[R(Q,\alpha)]^{SL(\alpha)}$
can be obtained using the following Corollary of the results from \cite{dw}:
\begin{proposition}\label{decomp_func}
 If $W\in V^{\perp}$ and $m_1S_1+m_2S_2+\cdots+m_tS_t$ 
is the sum of Jordan-H{\"o}lder factors of $W$ in $V^{\perp}$,
then $c_W = c_{S_1}^{m_1}c_{S_2}^{m_1}\cdots c_{S_t}^{m_t}$. 
\end{proposition}
\begin{proof}
We have a filtration 
$0=W_0\subseteq W_1\subseteq \cdots\subseteq W_d = W$ such that
$W_j\in V^{\perp}$, $W_j/W_{j-1}\cong S_p$,
$j=1,\cdots,d$. Applying \cite[Lemma~1]{dw}, we decompose $c_W$.
\end{proof}

Denote by $n_o$ the number of $c$-orbits in $I$;
then $n_o=2$ if $\Gamma=\widetilde{A_n}$ and $n_o=3$,
otherwise. For each orbit ${\mathcal O}_i,i=1,\cdots,n_o$,
denote by $P_i$ the product of $c_{E_i}$ over
the orbit; clearly this semi-invariant is of weight $\sigma$.

\begin{theorem}\label{homog_gen}\quad

{\bf 1}.The algebra 
${\bf k}[R(Q,p\delta)]^{SL(p\delta)}$
is generated by $c_{E_i}$, $i\in I$ and
$f_0,\cdots,f_p$.

{\bf 2}. A minimal system of generators of ${\bf k}[R(Q,p\delta)]^{SL(p\delta)}$
consists of $c_{E_i}$, $i\in I$ and $max(p+1-n_o,0)$ elements
from $f_0,\cdots,f_p$. If $p+1 \geq n_o$, then
these generators are algebraically independent;
otherwise, if $p=1,n_o=3$, then the generators fulfill 
a syzygy $c_1P_1+c_2P_2+c_3P_3=0,c_1,c_2,c_3\in{\bf k}^*$ 
and the ideal of syzygies is generated by this one. 
\end{theorem}

\begin{remark}
This statement is the same as \cite[Theorem~2.4]{skw}.
\end{remark}

\begin{proof}
Denote by $E_{\lambda},\lambda\in \Lambda\subseteq {\bf k}$ 
a 1-parameter family
of pairwise non-isomorphic simple homogeneous regular representations
of $Q$ of dimension $\delta$.
A generic representation of dimension $p\delta$ is locally
semi-simple and is isomorphic to $E_{\lambda_1}+\cdots+E_{\lambda_p}$.
By \cite{dw}, ${\bf k}[R(Q,\alpha)]^{SL(\alpha)}$ is generated by 
the semi-invariants
$c_W$ such that $W\in (E_{\lambda_1}+\cdots+E_{\lambda_p})^{\perp}$ 
for some collection $\lambda_1,\cdots,\lambda_p$; 
moreover, by Proposition \ref{decomp_func}
we can assume $W$ to be a simple object of this category, hence,
a $\sigma$-stable representation, by Proposition \ref{simple_is_stable}.
The $\sigma$-stable representations are the simple regular ones.
Note also that for $W$ being homogeneous simple, 
$c_W\subseteq 
{\bf k}[R(Q,p\delta)]^{(GL(p\delta))}_{\sigma}=\langle f_0,\cdots,f_p\rangle$.
So the assertion {\bf 1} is proved.

The generic stabilizer of $GL(p\delta)$ is $({\bf k}^*)^p$,
hence the generic stabilizer of $SL(p\delta)$ is $({\bf k}^*)^{p-1}$.
So $\dim{\bf k}[R(Q,p\delta)]^{SL(p\delta)}=
\dim R(Q,p\delta) - \dim SL(p\delta) + (p-1) = 
q_Q(p\delta) + \vert Q_0\vert +(p-1)=
n+p$, where $n=\vert Q_0\vert-1$.

Consider the semi-invariants $P_j,j=1,\cdots,n_o$.
Set $W_j=\sum_{i\in{\mathcal O}_j}E_i$.
Since ${\rm Hom}(W_j,E_i)$ is non-trivial if and only if
$i\in {\mathcal O}_j$, $P_j$ vanishes on $W_k$ if and only if
$k=j$. Hence, $P_i$ and $P_j$ are non-proportional for $i\neq j$.
Moreover, if $n_O=3$ and $p\geq 2$, then
$P_1,P_2,P_3$ are linearly independent because of the 
values of these on representations $W_1+W_2,W_1+W_3,W_2+W_3$.

Therefore, if $p+1 \geq n_o$, then the semi-invariants $P_1,\cdots,P_{n_o}$
are linear independent and 
${\bf k}[R(Q,p\delta)]^{SL(p\delta)}$ is generated by
$c_{E_i},i\in I$ and  $p+1-n_o$ elements of $f_0,\cdots,f_p$. 
One can see that $I$ consists of $n+n_O-1$ elements. 
So this system of generators
consists of $(n+n_O-1)+(p+1-n_O)=n+p$ elements,
hence, this is a minimal system of algebraically independent
generators.

Finally, if $p=1$, $n_O=3$, then $P_1,P_2$, and $P_3$
are non-proportional elements of two-dimensional vector space
${\bf k}[R(Q,p\delta)]^{(GL(p\delta))}_{\sigma}$, hence we get
a syzygy as in {\bf 2}. Since the number of generators
is $n+2=n+p+1$, and because our syzygy is of degree
1 by each of the generators, the assertion {\bf 2} is proved. 
\end{proof}

Now consider the general case: $\alpha = p\delta+\sum_{i\in I}p_ie_i$.
The quiver $E_Q$ is the union of 2 or 3 circular quivers; 
for $i\in I$ define by $n(i)$ and $p(i)$ the next and the
previous vertex of $E_Q$ so that $c(e_i)=e_{n(i)}, n(p(i))=i$. 
A subset $[k,l]=\{k,n(k),\cdots,l\}\subseteq I$ will be called an {\it arc}.
By Proposition \ref{D_V_prop}.1, each arc $[k,l]$ with $k\neq l$
yields a real Schur root $e_{k,l}=e_k+\cdots+e_l$;
pick a Schurian representation $E_{k,l}\in R(Q,e_{k,l})$. 
By Proposition 
\ref{slice_decomp}.3 and Algorithm \ref{alg_gen_lss_A_n}
the generic locally semi-simple decomposition of $\alpha$ is
$\alpha = p\delta+\sum_{[k,l]\in \Omega} m_{k,l}e_{k,l}$,
where $\Omega$ is a set of arcs such that for
different arcs $[k_1,l_1],[k_2,l_2]\in \Omega$

$\bullet$ either $[k_1,l_1]\cap[k_2,l_2]=\emptyset$ or
$[k_1,l_1]\subseteq [n(k_2),p(l_2)]$ or else $[k_2,l_2]\subseteq [n(k_1),p(l_1)]$

$\bullet$ if $m_{k_1,l_1} = m_{k_2,l_2}$, then $p_{k_1}\neq p_{k_2}$.

\begin{proposition}\label{dim_gen}
$\dim {\bf k}[R(Q,\alpha)]^{SL(\alpha)} = 
\dim {\bf k}[R(Q,p\delta)]^{SL(p\delta)} - \vert\Omega\vert$.
\end{proposition}

\begin{proof}
By Theorem \ref{lr_quiver}.3, ${\bf k}[R(Q,\alpha)]^{SL(\alpha)}\cong
{\bf k}[R(Q,p\delta)\oplus\bigoplus_{[k,l]\in\Omega}R(Q,e_{k,l})]^G$,
where $G\subseteq GL(p\delta)\times\prod_{[k,l]\in\Omega}GL(e_{k,l})$
consists of the elements with the product of determinants at any vertex 
being 1; moreover, generic $G$-orbit is closed.
Next, since $e_{k,l}$ are real Shur roots,
$G$ acts on $\bigoplus_{[k,l]\in\Omega}R(Q,e_{k,l})$
with an open orbit $G/K$, where $K$ is the kernel of that action.
Consequently, generic $K$-orbit is closed in $R(Q,p\delta)$
and $\dim {\bf k}[R(Q,\alpha)]^{SL(\alpha)}=
\dim {\bf k}[R(Q,p\delta)]^K$. Observe that $K$ acts on $R(Q,p\delta)$
as a subgroup $TSL(p\delta)\subseteq GL(p\delta)$, where 
$T$ is a central torus in $GL(p\delta)$ of dimension $\vert\Omega\vert$. 
One can easily check that $T$ acts effectively on $R(Q,p\delta)/\!\!/SL(p\delta)$,
hence $\dim {\bf k}[R(Q,p\delta)]^K = 
       \dim {\bf k}[R(Q,p\delta)/\!\!/SL(p\delta)]^T=
       \dim {\bf k}[R(Q,p\delta)/\!\!/SL(p\delta)] - \dim T = 
       \dim {\bf k}[R(Q,p\delta)]^{SL(p\delta)}  -\vert\Omega\vert$.
\end{proof}
 
By Theorem \ref{lr_quiver}.4,
${\bf k}[R(Q,\alpha)]^{SL(\alpha)}$ is isomorphic to 
$\bigoplus_{\chi\in \Lambda}
{\bf k}[R(Q,p\delta)]^{(GL(p\delta))}_{\chi}$,
where $\Lambda\subseteq{\bf Z}_+^{Q_0}$ consists of weights
such that ${\bf k}[R(Q,p\delta)]^{(GL(p\delta))}_{\chi}\neq 0$ and
for each $[k,l]\in\Omega$,
${\bf k}[R(Q,e_{k,l})]^{(GL(e_{k,l}))}_{\chi}\neq 0$.
By Theorem \ref{homog_gen}, ${\bf k}[R(Q,p\delta)]^{(GL(p\delta))}_{\chi}\neq 0$
implies $\chi=-\langle \;, e\rangle$, where 
$e\in \langle e_i,i\in I\rangle_{{\bf Z}_+}$.
So in order to determine $\Lambda$, we need to find the dimensions
$e\in \langle e_i,i\in I\rangle_{{\bf Z}_+}$ such that 
$E_{k,l}^{\perp}\cap R(Q,e)\neq\emptyset$ for any $[k,l]\in\Omega$.

\begin{proposition}\label{single}
$E_{k,l}^{\perp}\cap R(Q,e)\neq\emptyset$ iff
$e\in \langle e_{k,n(l)},e_i\vert i\in I,i\neq k,n(l)\rangle_{{\bf Z}_+}$.
\end{proposition}

\begin{proof}
Clearly, a necessary condition for $E_{k,l}^{\perp}\cap R(Q,e)\neq\emptyset$ 
is $\langle e_{k,l},e\rangle=0$. By Proposition\ref{isometry}
we have: 
\begin{equation}\label{Euler_induced}
\langle e_{k,l},\sum_{i\in I}q_ie_i\rangle = q_k - q_{n(l)}.
\end{equation}
Hence, the semi-group $\{e=\sum_{i\in I}q_ie_i,q_i\in{\bf Z}_+\vert 
\langle e_{k,l},e\rangle=0\}$
is generated by dimensions $e_{k,n(l)},e_i,i\in I\setminus\{k,n(l)\}$. 
So it is sufficient to check either of the equivalent
conditions ${\rm Hom}(E_{k,l},E)=0$ or  ${\rm Ext}(E_{k,l},E)=0$
for $E=E_{k,n(l)},E_i,i\neq k,n(l)$. For all $E$ with except of 
$E_{k,n(l)},E_l$ we have by (\ref{A_n_hom_hom}) and 
Proposition \ref{D_V_prop}.2:  ${\rm Hom}(E_{k,l},E)=0={\rm Hom}(E,E_{k,l})$.
On the other hand, for $E=E_{k,n(l)},E_l$,  (\ref{A_n_ext_ext})
and Proposition \ref{D_V_prop}.2 yield: ${\rm Ext}(E_{k,l},E)=0={\rm Ext}(E,E_{k,l})$.
\end{proof}

Let $J\subseteq I$ consist of elements being $k$ or $n(l)$
for an arc $[k,l]\in\Omega$. It can happen that
$J$ consists of less than $2\vert\Omega\vert$ elements because
there can be arcs like $[k,l]$ and $[n(l),m]$ in $\Omega$ such 
that their union is again an arc. So we can introduce
a new set $\Delta$ of arcs such that each arc from
$\Delta$ is a disjoint union of arcs from $\Omega$,
each arc from $\Omega$ is contained in an arc from $\Delta$,
and for any $[k_1,l_1],[k_2,l_2]\in\Delta$ we have:
$p(k_1)\neq l_2,n(l_1)\neq k_2$.

\begin{proposition}\label{lambda}
$\Lambda$ is generated by $\vert I\vert-\vert\Omega\vert$ elements
$\chi=-\langle \;,e\rangle$, where 
$e\in \{ e_i, e_{k,n(l)}\vert i\in I\setminus J,[k,l]\in\Delta\}$.
\end{proposition}

\begin{proof}
By formula (\ref{Euler_induced}) a necessary condition for a character 
$\chi=-\langle \;,\sum_{i\in I}q_ie_i\rangle$
to be in $\Lambda$ is $q_{n(l)}=q_k$ for any arc $[k,l]\in\Omega$.
Hence, the semigroup of dimension vectors meeting this
condition is generated by $e_i,i\in I\setminus J$ and
$e_{k,n(l)},[k,l]\in\Delta$. On the other hand, by Proposition \ref{single}
for each $e$ of this generators and for each arc $[k,l]\in\Omega$
there is a representation of dimension $e$ perpendicular to $E_{k,l}$.
It remains to note: 
$\vert J\vert = \vert \Delta\vert+\vert\Omega\vert$. 
\end{proof}

\begin{theorem}\label{main}
Let $\alpha=p\delta+\sum_{i\in I}p_ie_i,p>0$.
If $p=1$, $I$ consists of 3 orbits, and for each orbit 
at least two coefficients $p_i$ vanish, then  ${\bf k}[R(Q,\alpha)]^{SL(\alpha)}$
is a hypersurface; in all other cases it is a polynomial algebra.
\end{theorem}

\begin{proof}
Clearly, $A=\bigoplus_{\chi\in \Lambda}
{\bf k}[R(Q,p\delta)]^{(GL(p\delta))}_{\chi}$ is generated 
as an algebra by the 
subspaces ${\bf k}[R(Q,p\delta)]^{(GL(p\delta))}_{\chi}$, where
$\chi$ is $\sigma$ or a generator of $\Lambda$ from \ref{lambda}.
Moreover, for each orbit ${\mathcal O}_j\subseteq I$, 
$\delta$ can be obtained as a non-negative linear 
combination of the generators of $\Lambda$ corresponding
to  ${\mathcal O}_j$, hence $P_j$ is a corresponding product 
of generators of $A$.   
So $A$ is generated by $c_{E_i},i\subseteq I\setminus J$,
$c_{E_k}\cdots c_{E_{n(l)}}, [k,l]\in\Delta$, and
$max(p+1-n_o,0)$ elements from $f_0,\cdots,f_p$.
The number of these generators of $A$ is less
than the number of generators of ${\bf k}[R(Q,p\delta)]^{SL(p\delta)}$ by 
$\vert\Omega\vert$, hence by Proposition \ref{dim_gen}
if ${\bf k}[R(Q,p\delta)]^{SL(p\delta)}$ is polynomial algebra, 
then $A$ is and if ${\bf k}[R(Q,p\delta)]^{SL(p\delta)}$ is a hypersurface,
then $A$ is generated by $\dim A+1$ elements.
Moreover, in the latter case the unique relation between the generators is
$c_1P_1+c_2P_2+c_3P_3=0$. If for some ${\mathcal O}_j$, 
$P_j=c_{E_k}\cdots c_{E_{n(l)}}, [k,l]\in\Delta$,
then the relation says that this generator is redundant, 
so $A$ is in fact a polynomial algebra.
This happens precisely when $n(n(l))=k$ or equivalently,
$p_{n(l)}=0$ and $p_i\neq 0$ for all other $i\in{\mathcal O}_j$.
This completes the proof.
\end{proof}



\begin{thebibliography}{20}
\bibitem[BGP]{bgp} I.N.Bernstein, I.M.Gelfand, and V.A.Ponomarev,
Coxeter fuctors and Gabriel's theorem, Uspechi Mat. nauk {\bf 28} (1973), 19-33
(English transl. Rus. Math. Surveys {\bf 28} (1977), 17-32.)
\bibitem[DW]{dw} H. Derksen and J. Weyman, Semi-invariants of quivers and saturation 
for Littlewood-Richardson coefficients, J. AMS, {\bf 13} (2000), 3, 467-479.
\bibitem[GV]{gv} V. Gatti and E. Viniberghi, Spinors of 13-dimensional space,
Adv. in Math. {\bf 30} (1978), 137-155.
\bibitem[HH]{hh} R. Howe and R. Huang, Projective invariants of four subspaces, 
Adv. in Math. {\bf 118} (1996), 295-336.
\bibitem[Kac]{kac} V. Kac, Infinite root systems, representations of graphs,
and Invariant theory, II, J. Algebra {\bf 78} (1982), 141-162.  
\bibitem[Ki]{ki} A. D. King, Moduli of representations of finite dimensional 
algebras, Quart. J. Math. Oxford (2), {\bf 45} (1994), 515-530.
\bibitem[Kr]{kr} H. Kraft, Geometrische methode in der Invariantentheorie,
Vieweg, Braunschweig, 1984.
\bibitem[LBP]{lbp} L. Le Bruyn and C. Procesi, Semisimple representations of quivers,
Transactions of AMS, {\bf 317} (1990), 2, 585-598.
\bibitem[Lu1]{lu1} D. Luna, Slices \'etales, Bull. Soc. Math. France 
{\bf 33} (1973), 81-105.
\bibitem[Lu2]{lu2} D. Luna, Adh\'erence d'orbites et invariants,
Inv. Math. {\bf 29} (1975), 231-238.
\bibitem[Ma]{ma} Y. Matsushima, 
Espaces homog\`enes de Stein des groupes de Lie complexes, Nagoya
Math. J. 16 (1960), 205-218.
\bibitem[MF]{mf} D. Mumford and J. Fogarty, 
{\it Geometric Invariant Theory}, Springer-Verlag, 2nd. ed, 1982. 
\bibitem[Po]{po} V. L. Popov, A criterion of stability of actions of 
semisimple groups on factorial varieties, Izv. Akad. Nauk USSR,
{\bf 34} (1970), 523-531.
\bibitem[Ri]{ri} C. M. Ringel, Rational invariants of the tame quivers, 
Inv.math. {\bf 58}, (1980), 217-239.
\bibitem[Sch]{sch} A. Schofield, Semi-invariants of quivers, J. London Math. Soc.
{\bf 43} (1991), 383-395.
\bibitem[Sh]{sh} D. Shmelkin, On spherical representations of quivers and 
generalized complexes, Transf. groups, {\bf 7 } (2002), 1, 87-106.
\bibitem[SK]{sk} M. Sato and T. Kimura, A classification of irreducible 
prehomogeneous vector spaces and their relative invariants, Nagoya J. Math. 
{\bf 65} (1977), 1-155.
\bibitem[SkW]{skw} A. Skowronsky and J. Weyman, 
The algebras of semi-invariants of quivers, Transformation groups {\bf 5} (2000), 4, 361-402.
\bibitem[SchW]{schw} G. W. Schwarz and D. L. Wehlau, Invariants of
four subspaces, Ann. Inst. Fourier, Grenoble, {\bf 48} (1998), 3, 667-697. 
\end{thebibliography}
\end{document}